\documentclass[12pt,fullpage]{article}

\usepackage{amssymb}
\usepackage{amsmath}
\usepackage{epsfig}
\usepackage{caption}
\usepackage{psfrag}
\usepackage{array}
\usepackage{latexsym}
\usepackage{amsbsy,amsfonts,amssymb,amscd}
\usepackage{graphicx, color}

\input epsf
\numberwithin{equation}{section}

\newtheorem{Lemma}{Lemma}[section]
\newtheorem{Cor}[Lemma]{Corollary}
\newtheorem{Th}[Lemma]{Theorem}

\newtheorem{Prop}[Lemma]{Proposition}

\newtheorem{Rk}[Lemma]{Remark}
\newtheorem{Conj}{Conjecture}

\newcommand{\be}{\begin{equation}}
\newcommand{\ee}{\end{equation}}
\newcommand{\baa}{\begin{array}}
\newcommand{\eaa}{\end{array}}
\newcommand{\ba}{\begin{eqnarray}}
\newcommand{\ea}{\end{eqnarray}}
\newcommand{\ds}{\displaystyle}
\newcommand{\wto}{\rightharpoonup}

\def\R{{\mathbb R}}
\def\N{{\mathbb N}}

\def\L{{\mathcal L}}
\def\epsilon{\varepsilon}
\def\trait (#1) (#2) (#3){\vrule width #1pt height #2pt depth #3pt}
\def\fin{\hfill\trait (0.1) (5) (0) \trait (5) (0.1) (0) \kern-5pt \trait (5) (5) (-4.9) \trait (0.1) (5) (0)}
\def\2{C^{1,2}(\R\times\R^N)}
\def\to{\rightarrow}

\def\tilde{\widetilde}
\def\.{\cdot}

\begin{document}

\date{}
\title{\bf{A variational approach to reaction diffusion equations with forced speed in dimension 1}}
\author{Juliette Bouhours$^{\hbox{\small{1}}}$, Gr\'egoire Nadin$^{\hbox{\small{1}}}$\\
\\
\footnotesize{$^{\hbox{1}}$Sorbonne Universit\'es, UPMC Univ Paris 06, UMR 7598, Laboratoire Jacques-Louis Lions, F-75005, Paris, France,}\\
\footnotesize{CNRS, UMR 7598, Laboratoire Jacques-Louis Lions, F-75005, Paris, France.}
}

\maketitle

\begin{abstract}
We investigate in this paper a scalar reaction diffusion equation with a nonlinear
reaction term depending on $x-ct$. Here, $c$ is a prescribed parameter
modelling the speed of climate change and we wonder whether a population
will survive or not, that is, we want to determine the large-time behaviour
of the associated solution.

\noindent This problem has been solved recently when the nonlinearity is of KPP
type. We consider in the present paper general reaction terms, that are
only assumed to be negative at infinity. Using a variational approach,
we construct two thresholds $0<\underline{c}\leq \overline{c} <\infty$
determining the existence and the non-existence of travelling waves.
Numerics support the conjecture $\underline{c}=\overline{c}$. We then
prove that any solution of the initial-value problem converges at large
times, either to $0$ or to a travelling wave. In the case of bistable
nonlinearities, where the steady state $0$ is assumed to be stable, our
results lead to constrasting phenomena with respect to the KPP framework.
Lastly, we illustrate our results and discuss several open questions
through numerics.
\end{abstract}


\vskip 20pt
\section*{Acknowledgments}
The research leading to these results has received funding from the European Research Council under the European Union's Seventh Framework Program (FP/2007-2013) / ERC Grant Agreement n.321186 - ReaDi- Reaction-Diffusion Equations, Propagation and Modelling. Juliette Bouhours is funded by a PhD fellowship "Bourse hors DIM" of the "R\'egion Ile de France".

\section{Introduction and main results}\label{intro}

\subsection{Motivation: models on climate change}

Reaction diffusion problems are often used to model the evolution of biological species. In 1937, Kolmogorov, Petrovskii and Piskunov in \cite{KPP}, Fisher in \cite{Fisher} used reaction diffusion to investigate the propagation of a favourable gene in a population. One of the main notions introduced in \cite{KPP, Fisher} is the notion of travelling waves, i.e solution of the form $\ds{u(t,x)=U(x-ct)}$ for $x\in\R$, $t>0$ and some constant $c\in\R$. 
Since then a lot of papers have been dedicated to reaction diffusion equations and travelling waves in settings modelling all sorts of phenomena in biology.\\ 
In this paper we are interested in the following problem,
\be\label{problem}\tag{P}
\begin{cases}
u_t-u_{xx}=f(x-ct,u), & x\in\R, t>0,\\
u(0,x)=u_0(x), &x\in\R,
\end{cases} 
\ee
where $u_0$ is bounded, nonnegative and compactly supported.\\
This problem has been proposed in \cite{BDNZ} to model the effect of climate change on biological species. 
In this setting $u$ is the density of a biological population that is sensitive to climate change. We assume that the North Pole is found at $+\infty$ 
whereas the equator is at $-\infty$, which gives a good framework to study the effect of global warming on the distribution of the population. 
The dependence on $z$ in the reaction term takes into account the notion of favourable/unfavourable area depending on the latitude for populations which are sensitive 
to the climate/temperature of the environment. The constant $c$ can be seen as the speed of the climate change. In such a setting, one will be interested to know 
when the population can keep track with its favourable environment despite the climate change and thus persists at large times. In \cite{BDNZ} Berestycki et al studied the 
existence of non trivial travelling wave solutions converging to $0$ at infinity in dimension 1 when $f$ satisfies the KPP property: $\ds{s\in\R^+\mapsto f(z,s)/s}$ 
is decreasing for all $z\in\R$. 
As they converge to $0$ as $z\to \pm \infty$, such homoclinic solutions are sometimes called travelling pulses in the literature, 
in contrast with heteroclinic travelling waves. 
We will use in this paper the same terminology as in \cite{BDNZ} and call travelling waves such homoclinic solutions (see equation (\ref{elliptproblem}) below
for a precise definition).
Berestycki et al proved that in this framework, the persistence of the population depends on the sign of the principal eigenvalue of the linearized 
equation around the trivial steady state 0. Their results have been extended by Berestycki and Rossi to $\R^N$ in \cite{BR1} and to infinite cylinders in \cite{BR2}. In \cite{Vo} Vo
studies the same type of problem with more general classes of unfavourable media toward infinity. 

A similar model was developed by Potapov and Lewis in \cite{PL} and by Berestycki, Desvillettes and Diekmann in \cite{BDD} in order to investigate a two-species competition system facing climate change. These papers studied the effect of the speed of the climate change on the coexistence between the competing species. In \cite{BDD} the authors pointed out 
the formation of 
a spatial gap between the two species when one is forced to move forward to keep up with the climate change and the other has limited invasion speed. 
The persistence of a specie facing climate change was also investigated mathematically through an integrodifference model by Zhou and Kot in \cite{ZK}.
 
The particularity of all these papers is the KPP assumption for the 
reaction term, where the linearized equation around 0 determines the behaviour of the solution of the nonlinear equation.
As far as we know, such questions were only investigated numerically for other types of nonlinearities by Roques et al in \cite{BRRK}, where the authors were mainly interested in the effects of the geometry of the domain (in dimension 2) on the persistence of the population considering KPP and bistable nonlinearities.


\subsection{Framework}
In this paper we are interested in this persistence question, when the evolution of the density of the population is modelled by a reaction diffusion equation, with more general hypotheses on the nonlinearity $f$. Indeed we point out that we consider general nonlinearities $f$, without assuming $f$ to satisfy the KPP property. \\
We will assume that $f$ is a Carath\'eodory function, i.e 
\begin{align*}
s\mapsto f(z,s) \text{ is continuous for all } z\in\R,\\
z\mapsto f(z,s) \text{ is measurable for all } s\in\R,
\end{align*}
satisfying the following hypotheses,
\be\label{f0}
f(z,0)=0,
\ee
\be\label{flip}
s\mapsto f(z,s) \text{ is Lipschitz-continuous and of class } \mathcal{C}^1\text{ uniformly with respect to } z\in\R,
\ee
\be\label{ubdd}
\exists M>0 \hspace{0.1cm},  f(z,s)\leq 0, \hspace{0.5cm} \forall s\in\R\backslash(0,M) \text{ and } z\in\R,
\ee
\be\label{favbdd}
\exists R>0, \delta>0, \hspace{0.2cm} f(z,s)\leq -\delta s, \hspace{0.5cm} \forall |z|>R, s\in(0,M).
\ee
\be \label{hypenergy}
\exists u\in H^1(\R), \hspace{0.2cm} E_0[u]:= \int_\R \left(\frac{u_z^2}{2}-F(z,u)\right)dz<0, \quad \text{with } F(z,s):= \int_0^sf(z,t)dt,
\ee
Assumption \eqref{f0} means than when the population vanishes then no reaction takes place, i.e 0 is a steady state of the problem which corresponds 
to the extinction of the population. Hypothesis \eqref{ubdd} models some overcrowding effect: the resources being limited, the environment becomes unfavourable when the population grows above some threshold $M>0$. Assumption \eqref{favbdd} gives information on 
the boundedness of the favourable environment and postulates that outside a bounded region the environment is strictly unfavourable.
Lastly, hypothesis \eqref{hypenergy} will imply that the equation with $c=0$ admits a non-trivial stationary solution which is more stable, in a sense, than 
the equilibrium $u\equiv 0$.
\bigskip\\
We will use the following weighted spaces throughout this paper:
$$L_c^2(\R):=L^2(\R,e^{cz}dz),\: H_c^1(\R):=H^1(\R,e^{cz}dz),\: H_c^2(\R):=H^2(\R,e^{cz}dz).$$
\subsection{Main results}
Up to a change of variable ($z:=x-ct$) Problem \eqref{problem} is equivalent to
\be\label{problembis}\tag{$\tilde{P}$}
\begin{cases}
u_t-u_{zz}-cu_z=f(z,u), & z\in\R, t>0,\\
u(0,z)=u_0(z), &z\in\R,
\end{cases} 
\ee
In our paper we investigate the existence of travelling wave solutions of \eqref{problem}, i.e nonnegative solution of the form $u(t,x)=U(x-ct)$ for all $x\in\R$, $t>0$ with $\ds{U\not\equiv 0}$, $\ds{U(\pm\infty)=0}$ . This particular solutions are non trivial solutions of the following stationary problem 
\be\label{elliptproblem} \tag{S}
\begin{cases}
-U_{zz}-cU_z=f(z,U),& z\in\R,\\
U(z)\geq 0, &z\in\R,\\
U(\pm\infty)=0.
\end{cases}
\ee
Solutions of \eqref{elliptproblem} are also the stationary solutions of Problem \eqref{problembis} and notice that 0 is a solution of \eqref{elliptproblem} but not a travelling wave solution. We have the following theorem,

\begin{Th}\label{notrivialTW}
Under hypotheses \eqref{f0}-\eqref{hypenergy}, there exist $\overline{c} \geq \underline{c}>0$, such that 
\begin{itemize}
\item for all $c\in(0,\underline{c})$, \eqref{problem} has a travelling wave solution $\ds{U_c\in H_c^1(\R)}$, 
\item For all $c>\overline{c}$, \eqref{problem} has no travelling wave solution, that is 0 is the only solution of \eqref{elliptproblem}.
\end{itemize}
\end{Th}
The proof of Theorem \ref{notrivialTW} is based on a variational approach, used in \cite{LMN} to prove the existence of travelling  front for gradient like systems of equations. We use the same variational formula but in the case of scalar equations and when $f$ depends on $z=x-ct$. Namely we introduce an energy $E_c$ defined by \eqref{Energy} such that for all $c>0$ if there exists $u\in H_c^1(\R)$ with $E_c[u]<0$, then there exists a travelling wave solution constructed as a minimiser of $E_c$. Theorem \ref{notrivialTW} is derived from a continuity argument using assumption \eqref{hypenergy}.
\bigskip\\
Then we will be interested in the convergence of the Cauchy problem.
 
\begin{Th}\label{uniqlimit}
Let $u_0\in H^2(\R)$ and $u_0$ non-negative and compactly supported. Then the unique solution $u$ of \eqref{problem} satisfies $u\in L^2([0,T[,H_c^1(\R))$, $u_t\in L^2([0,T[,L_c^2(\R))$, for all $T>0$, and $\ds{t\mapsto u(t,\cdot-ct)}$ converges to a solution of \eqref{elliptproblem} in $H_c^2(\R)$ as $t\to+\infty$.
\end{Th} 
Note that the limit of $u$ in the previous theorem could be either the trivial solution 0 or a travelling wave. And if $\ds{t\mapsto u(t,\cdot-ct)}$ converges to 0 as $t\to+\infty$ this implies that the population goes exctinct, whereas if $\ds{t\mapsto u(t,\cdot-ct)}$ converges to $U_c>0$ non trivial solution of \eqref{elliptproblem} as $t\to+\infty$ this means that there is persistence of the population and convergence to a travelling wave solution.
\bigskip\\ 
After proving these two main theorems, we study the existence of travelling wave solutions and the behaviour of the solution of the Cauchy problem \eqref{problembis} depending on the linear stability of 0. Then we study the solution $u$ of \eqref{problembis} for particular $f$, $\delta$ and $c$.
\begin{itemize}
\item We prove that, as in the KPP framework, when 0 is linearly unstable the solution $u$ of \eqref{problembis} converges to a travelling  
wave solution. We also exhibit some bistable-like frameworks, that is, frameworks where 0 is linearly stable,
admitting travelling wave solutions (see Section \ref{sec:bistable}) but for which solutions could converge to 0 or a travelling wave depending on the initial datum.
This emphasises the particularity of the KPP framework where the linearized equation near $u\equiv 0$ determines the existence of travelling wave solutions, which is not true in the general framework.
\item In the last section we first study numerically the existence of a speed threshold such that the population survives if $c$ is below this threshold, 
i.e the solution of the Cauchy problem \eqref{problembis} converges toward travelling waves for large times, while the population dies if the speed $c$ is above this threshold, 
i.e the solution of the Cauchy problem \eqref{problembis} converges to 0 for large times, for $f$ KPP, monostable and bistable in the favourable area. For such nonlinearities, we believe that we chose our initial data in such a way that they lie in the basin of attraction of a travelling wave solution, when it exists. Thus in view of the numerical results, we state the following conjecture:
\begin{Conj}
Let $\underline{c}\leq\overline{c}$ be defined by Theorem \ref{notrivialTW} then $\ds{\underline{c}=\overline{c}}$.
\end{Conj}
We also plot the shape of the profile for different values of the parameter $\delta$ and $f$ bistable.
\bigskip\\ 
Then we give an example of nonlinearity $f$ such that there exist several locally stable travelling wave solutions in Proposition \ref{energyneg} and illustrate this result with numerical simulations displaying the shape of the profile for different initial conditions. 
\end{itemize}
\textbf{Organisation of the paper}
\bigskip\\
Theorem \ref{notrivialTW} concerning the stationary framework is proved using a variational method in section \ref{stationnary}. Sections \ref{CauchyPb} and \ref{stab0} are devoted to the study of the Cauchy problem \eqref{problem}. We prove Theorem \ref{uniqlimit} in section \ref{CauchyPb} and discuss the linear stability of 0 and its consequences on the convergence of the Cauchy problem in section \ref{stab0}. We give some examples and discuss possible improvements of our results with numerical insight in section \ref{examples}.


\section{A variational approach to travelling waves}\label{stationnary}

The variational structure of travelling wave solutions of homogeneous
reaction-diffusion equations is known since the pioneering work of Fife and McLeod \cite{FM}.
However, this structure has only been fully exploited quite recently in order to
derive existence and stability results for travelling waves in bistable equations in
parallel by Heinze \cite{Heinze}, Lucia, Muratov and Novaga \cite{M04, LMN04, LMN} and then by Risler \cite{Risler} for gradient systems (see
also \cite{GR, GJ} for various other applications).
The situation we consider in the present paper is different. First, we deal with
heterogeneous reaction-diffusion equations. The homogeneity was indeed a difficulty
in earlier works, since the invariance by translation caused a lack of
compactness. Here, the behaviour of the nonlinearity at infinity will somehow trap
minimising sequences in the favourable habitat where $ f $ is positive. Second, we
consider general nonlinearities, including monostable ones. The variational approach
is not a relevant tool in order to investigate such equations when the coefficients
are homogeneous since travelling waves do not decrease sufficiently fast at infinity
and thus have an infinite energy. Here, again, the behaviour of the nonlinearity at
infinity forces an admissible exponential decay and we could thus define an energy and make use of it.
\bigskip\\
We are interested in the existence of travelling wave solution of equation \eqref{problem}, i.e $\ds{u(t,x)=U(x-ct)=U(z)}$,  and $U$ is a solution of the ordinary differential equation
\begin{equation*}\begin{cases}
-U''-cU'=f(z,U), & z\in\R,\\
U>0 &\text{in }\R,\\
U(z)\to0 &\text{as } |z|\to+\infty.
\end{cases}\end{equation*}
\bigskip\\
To study existence of non trivial travelling waves, we introduce the energy functional defined as follow
\be\label{Energy}
E_c[u]=\int_\R e^{cz}\left\{ \frac{u_z^2}{2}-F(z,u)\right\}dz, \hspace{0.5cm} \forall u\in H_c^1(\R),
\ee
where $H_c^1(\R)=H^1(\R,e^{cx}dx)$ and 
$$F(z,s)=\int_0^sf(z,t)dt.$$ 
One can notice that \eqref{f0} and \eqref{flip} ensure that $\ds{\int_\R F(z,u)e^{cz}dz}$ is well defined for all $u\in H_c^1(\R)$. 
\vspace{0.5cm}
\noindent We start by proving the first part of the Theorem and by pointing out the link between solutions of \eqref{elliptproblem} and the functional $E_c$.
\begin{Lemma}\label{critpointsol}
Consider $u\in H^1_c(\R)$ a nonnegative function. Then $u$ is a critical point of the energy functional $E_c$ if and only if $u$ is a solution of \eqref{elliptproblem}. Moreover $u\in W^{2,p}_{\text{loc}}(\R)$, for all $p>1$.
\end{Lemma}
{\bf Proof}: The first part of the proof is classical. Let $u$ be a critical point of $E_c$. Standard arguments yield that $E_c$ is G\^ateaux-differentiable and that its differential at $u$ is given, for all $w\in H_{c}^1(\R)$, by
\be\label{diffenergy}
dE_c[u](w)=\int_\R e^{cz}\left\{u_zw_z-f(z,u)w\right\}dz.
\ee
This quantity is clearly continuous with respect to $u\in H^1_c (\R)$ as a linear form over $H^1_c(\R)$ since $f$ is Lipschitz-continuous by hypothesis (\ref{flip}).
Moreover letting $v(z):=u(z)e^\frac{cz}{2}$ for all $z\in\R$ and $w\in C^\infty$ compactly supported, then
$$v''=\frac{c^2}{4}v-f(z,e^{-\frac{cz}{2}}v)e^{\frac{cz}{2}},\quad \text{in }\R$$
and 
$$\frac{v''(z)}{v(z)}\geq \delta+\frac{c^2}{4},\quad \text{if } z\leq-R.$$
As $u\in H_c^1(\R)$, $v(z)\to0$ as $z\to-\infty$ and we can apply \cite[Lemma 2.2]{BR1} we have that
$$v(z)e^{-\sqrt{\frac{c^2}{4}+\delta}z}\underset{z\to-\infty}{\to}0,$$
which implies that 
$$u(z)\leq e^{\gamma z},\quad \forall z\leq R^-,$$
for some $R^-<-R$ and $\gamma>0$. This implies that $u(z)\to0$ as $z\to-\infty$ and as $u\in H_c^1(\R)$, $u(z)\to0$ as $z\to+\infty$ and $u$ is a solution of \eqref{elliptproblem}. \\
Now let $u\in H_c^1(\R)$ be a solution of \eqref{elliptproblem} then $u$ is a critical point of $E_c$, as $u(\pm\infty)=0$.
Thus $u\in H_c^1(\R)$ is a critical point of $E_c$ iff $u$ is a weak solution of \eqref{elliptproblem}.\\
Moreover $u$ is in $H^1_{loc} (\R)$, and thus continuous, and so is $v$. It follows from the equation on $v$ that it is $\mathcal{C}^2$, from which the conclusion follows.
 \fin
\bigskip\\
Let us state a Poincar\'e type inequality that will be useful in the sequel, which is due to \cite{LMN}.
\begin{Lemma}[Lemma 2.1 in \cite{LMN}]\label{lemmapoincareH1c}
For all $u\in H_c^1(\R)$,
\be\label{poincareH1c}
\frac{c^2}{4}\int_\R e^{cz}u^2dz\leq \int_\R e^{cz}u_z^2dz.
\ee
\end{Lemma}
Now notice that we can always assume that a global minimiser of $E_c$ is bounded and nonnegative.

\begin{Rk}\label{Ubdd} 
Considering $\ds{\tilde{u}=\min\left\{u,M\right\}}$, we have
\begin{align*}
F(z,u)-F(z,\tilde{u})=\int_{\tilde{u}}^uf(z,s)ds&=\begin{cases} 0 &\text{if } u<M,\\
							\ds{\int_M^uf(z,s)ds} &\text{otherwise},
					\end{cases}\\
				&\leq 0. \hspace{1cm}
\end{align*}
Thus
\begin{align*}
E_c[u]\geq \int_\R e^{cz}\left\{ \frac{u_z^2}{2}-F(z,\tilde{u})\right\}dz\geq \int_\R e^{cz}\left\{ \frac{\tilde{u}_z^2}{2}-F(z,\tilde{u})\right\}dz=E_c[\tilde{u}].
\end{align*}
As we want to minimise the energy functional, $\tilde{u}$ will always be a better candidate than $u$.

\noindent Similarly, taking $\tilde{u} = \max \{ 0, u\}$ instead of $u$ gives a lower energy.
\end{Rk}
Hypothesis \eqref{f0} ensures that $E_c(0)=0$ and thus $\underset{u\in H_c^1(\R)}{\inf}E_c[u]\leq 0$. Moreover, the following lemma yields that $\underset{u\in H_c^1(\R)}{\inf}E_c[u]>-\infty$.

\begin{Lemma}\label{Eminore}
For all $c>0$, there exists $C>0$ such that for all $u\in H_c^1(\R)$, $E_c[u]\geq -C$.
\end{Lemma}

\noindent {\bf Proof}: We can assume that $0\leq u\leq M$ using Remark \ref{Ubdd}. For all $u\in H_c^1(\R)$, using assumption \eqref{favbdd},

\begin{align*}
E_c[u] \geq \int_{-R}^Re^{cz}\left\{\frac{u_z^2}{2}-F(z,u)\right\}dz+\int_{\R\backslash (-R,R)}e^{cz}\left\{\frac{u_z^2}{2}+\frac{\delta u^2}{2}\right\}dz\geq E^R_c[u],
\end{align*}
where $\ds{E_c^R[u]=\int_{-R}^Re^{cz}\left\{\frac{u_z^2}{2}-F(z,u)\right\}dz}$. This implies that 
$\ds{\underset{u\in H_c^1(\R)}{\inf }E_c[u] \geq \underset{u\in H_c^1(-R,R)}{\inf }E_c^R[u]}$. 
Using the assumptions on $f$, there exists $C_0>0$ such that $-F(z,s)>-C_0$ for all $z\in (-R,R)$ and $s\in[0,M]$. Thus there exists $C>0$ such that 
${E_c^R[u]\geq -C}$ for all $u\in H^1(-R,R)$ and then 
\begin{equation*}
\underset{u\in H_c^1(\R)}{\inf }E_c[u] \geq -C.
\end{equation*}
 \fin

\begin{Prop}\label{minreached}
There exists $u_\infty\in H_c^1(\R)$ such that $\ds{E_c[u_\infty]=\underset{u\in H_c^1(\R)}{\min}E_c[u]}$.
\end{Prop}
To prove Proposition \ref{minreached} we consider $(u_n)_n$ a minimising sequence of $E_c$ in $H_c^1(\R)$, i.e such that 
$\ds{E_c[u_n]\to \underset{u\in H_c^1(\R)}{\inf} E_c[u]>-\infty}$ as $n\to+\infty$. In view of Remark \ref{Ubdd} we can assume that $u_n$ 
is bounded for $n$ large enough.

\begin{Lemma}\label{unbdd}
There exist $N\in\N$, $C_1>0$, locally bounded with respect to $c$, such that for all $n>N$, 
$$\ds{\lVert u_n\rVert^2_{H_c^1(\R)}=\int_\R e^{cz}\left\{(u_n)_z^2+u_n^2\right\}dz \leq \frac{1+C_1}{\min\{\frac{1}{2},\frac{\delta}{2}\}}}.$$
\end{Lemma}

\noindent {\bf Proof of Lemma \ref{unbdd}}: For all $u\in H_c^1(\R)$, bounded,
\begin{align*}
E_c[u]&\geq \int_{-R}^Re^{cz}\left\{\frac{u_z^2}{2}-F(z,u)\right\}dz+\int_{\R\backslash (-R,R)}e^{cz}\left\{\frac{u_z^2}{2}+\frac{\delta u^2}{2}\right\}dz\\
	&\hspace{0.5cm} = \int_{-R}^Re^{cz}\left\{-F(z,u)-\frac{\delta u^2}{2}\right\}dz+\int_{\R}e^{cz}\left\{\frac{u_z^2}{2}+\frac{\delta u^2}{2}\right\}dz\\
	&\hspace{0.5cm}\geq -C_1+\min\{\frac{1}{2},\frac{\delta}{2}\}\lVert u\rVert^2_{H_c^1(\R)},
\end{align*}
where $C_1 = -\frac{1}{c}(C_0+\frac{\delta M^2}{2}) (e^{cR}-e^{-cR})$, with $C_0$ as in the proof of Lemma \ref{Eminore}.
Moreover as $\ds{E_c[u_n]\to \underset{u\in H_c^1(\R)}{\inf} E_c[u]\leq 0 = E_c [0]}$, there exists $N\in\N$ such that for all $n>N$, ${E_c[u_n]\leq 1}$. 
Then using the previous computation we obtain the Lemma. \fin
\bigskip\\
One can now prove Proposition \ref{minreached}.
\bigskip\\
\noindent {\bf Proof of Proposition \ref{minreached}}: From Lemma \ref{unbdd}, if $(u_n)$ is a minimising sequence of $E_c$ in $H_c^1(\R)$ then  
$(u_n)$ is bounded in $H^1_c(\R)$. Thus up to a subsequence $(u_n)$ converges weakly to some $u_\infty\in H_c^1(\R)$. One has:
\be\label{ineqnormeH1c}
\int_\R e^{cz}(u_\infty)_z^2dz\leq \underset{n\to+\infty}{\liminf}\int_\R e^{cz}(u_n)_z^2dz.
\ee
Moreover as $u_n\in(0,M)$, classical Sobolev injections yield that
\begin{gather}
u_n\to u_\infty \textrm{ in } C_{loc}(\R) \textrm{ as } n\to+\infty \label{convvnC}.
\end{gather}
As $F$ is bounded, the dominated convergence theorem gives, for all $T\in\R$,
\be\label{dominatedconvF}
\int_{-\infty}^Te^{cz}F(z,u_n)dz\to \int_{-\infty}^T e^{cz}F(z,u_\infty)dz \quad \text{as } n\to+\infty.
\ee
Thus, as $-\int_T^{+\infty} e^{cz}F(z,u_n)dz\geq0$, for all $T>R$,
\begin{align*}
\underset{n\to+\infty}{\liminf}E_c[u_n]&\geq \int_\R e^{cz}\frac{(u_\infty)_z^2}{2}dz+\int_{-\infty}^T-e^{cz}F(z,u_\infty)dz\\
							&=E_c[u_\infty]+\int_T^{+\infty}e^{cz}F(z,u_\infty)dz\\
							&\geq E_c[u_\infty]-\int_T^{+\infty} C e^{cz} u_\infty^2dz,
\end{align*}
the last inequality following from \eqref{flip}. As, for all $\epsilon>0$, there exists $T>R$ such that $\int_T^{+\infty} C e^{cz} u_\infty^2dz<\epsilon$, 
since $u_\infty\in H_c^1$, we have 
$$E_c[u_\infty]\leq \underset{n\to+\infty}{\liminf} E_c[u_n]=\underset{u\in H_c^1(\R)}{\inf}E_c[u],$$
and the Proposition is proved. \fin
\bigskip\\
We have proved that the minimum is reached in $H_c^1$. This implies that there exists a solution $U$ of \eqref{elliptproblem} 
such that $E_c[U]=\underset{u\in H_c^1}{\inf} E_c[u]$. 
\bigskip\\

\begin{Prop}\label{continc}
The function $c\geq0\mapsto  \underset{u\in H^1_c(\R)}{\inf}E_c[u]$ is continuous.
\end{Prop}

\noindent {\bf Proof}: Let $(c_n)_n$ be a sequence in $\R$ such that $c_n\to c$ as $n\to+\infty$. From Proposition \ref{minreached} we know that for all $n\in\N$ there exists $u_n\in H^1_{c_n}(\R)$ 
such that  $\ds{\underset{u\in H^1_{c_n}(\R)}{\inf}E_{c_n}[u]=E_{c_n}[u_n]}$. Let $v_n:= e^\frac{c_n z}{2} u_n \in H^1(\R)$ and notice that 
$$E_{c_n}[u_n]= \int_\R \Big\{\frac{(v_n)_z^2}{2}+\frac{c_n^2}{8}v_n^2-e^{c_n z}F(z,e^{-\frac{c_nz}{2}}v_n) \Big\} dz=:\tilde{E}_{c_n}[v_n].$$ 
Moreover the sequence $(v_n)_n$ is uniformly bounded in $H^1(\R)$ by Lemma \ref{unbdd}, as $(c_n)_n$ is uniformly bounded, thus up to a subsequence, $v_n\wto v_\infty$ 
weakly in $H^1(\R)$ as $n\to +\infty$. Moreover, as $f$ is Lipschitz-continuous and $f(z,0)=0$ for all $z\in\R$, there exists a constant $C>0$ such that 
$|f(z,s)|\leq C s$ for all $s\geq 0$ and $z\in\R$. Integrating this inequality, one gets 
$e^{c_n z}| F(z,e^{-\frac{c_n z}{2}} v) |\leq C \frac{v^2}{2}$ for all $z\in\R$, $v\in H^1(\R)$ and $n\in\mathbb{N}$. It follows from the dominated convergence theorem and the continuity of $f$ in $s$ that 
$\lim_{n\to +\infty} \int_\R e^{c_n z} F(z,e^{-\frac{c_n z}{2}} v)dz = \int_\R e^{c z} F(z,e^{-\frac{c z}{2}} v)dz$ for all $v\in H^1(\R)$. Hence, $\tilde{E}_{c_n}[v]\to\tilde{E}_c[v]$ as $n\to+\infty$. As $v_n$ is a minimiser, for all $v\in H^1(\R)$,
\be
\tilde{E}_{c_n}[v_n]\leq \tilde{E}_{c_n}[v].
\ee
Passing to the limit and using the same arguments as in the proof of Proposition \ref{minreached} we obtain
\be
\tilde{E}_c[v_\infty]\leq \underset{n\to+\infty}{\liminf} \tilde{E}_{c_n}[v_n]\leq \tilde{E}_c[v], \quad \forall v\in H^1(\R).
\ee
This implies that $\ds{\tilde{E}_c[v_\infty]= \underset{v\in H^1(\R)}{\inf} \tilde{E}_c[v]}$, and letting $u_\infty=e^{-\frac{cz}{2}}v_\infty$ we get the Proposition. \fin
\vspace{0.2cm}\\
By the continuity of $c\mapsto  \underset{u\in H^1_c(\R)}{\inf}E_c[u]$ and using Proposition \ref{minreached}, the first part of Theorem \ref{notrivialTW} is proved.
\bigskip\\
This proposition will prove the second part of the Theorem
\begin{Prop}\label{cbar}
There exists $\overline{c}>0$ such that for all $c>\overline{c}$, 0 is the only solution of equation \eqref{elliptproblem}.
\end{Prop}
\textbf{Proof: }Define 
\be\label{gKPP}
g(z,u)=\left(\underset{s\geq u}{\sup}\frac{f(z,s)}{s}\right)\times u, \hspace{0.5cm} \forall z\in\R,u\in\R^+.
\ee
Then $g$ satisfies the following assumptions:
\begin{gather}
g(z,0)=0, \hspace{0.5cm} \forall z\in\R,\label{g0}\\
u\mapsto g(z,u) \text{ is Lipschitz-continuous uniformly with respect to } z\in\R, \label{glip}\\
g(z,s)\leq0, \hspace{0.5cm} \forall z\in\R, s\geq M,\label{gsolbdd}\\
\forall z\in\R,\quad u\mapsto \frac{g(z,u)}{u} \text{ is decreasing },\label{gKPP}\\
g(z,u)\leq-\delta u, \hspace{0.5cm} \forall |z|>R.\label{gfavbdd}
\end{gather}
Hence we know from \cite[Theorem 3.2]{BR1} that there exists $\overline{c}>0$ such that if $v$ is a solution of 
\begin{equation}\label{assKPP}
-v_{zz}-cv_z=g(z,v) \quad \text{in } \R,
\end{equation}
for $c>\overline{c}$, then $v\equiv 0$.
Moreover for all $z\in\R$, $s\in\R$, $g(z,s)\geq f(z,s)$. Take $c>\overline{c}$ and let $u$ a solution of \eqref{elliptproblem}, then $u$ is a subsolution of the associated KPP equation, i.e
$$-u''-cu'\leq g(z,u) \quad \text{in } \R.$$
Let $M>0$ be as in condition \eqref{ubdd}, then $w(z)=M$ for all $z\in\R$ is a super solution of the associated KPP problem, i.e
$$-w''-cw'\geq g(z,w) \quad \text{in } \R,$$
and we can take $M$ large enough such that $u\leq M$ in $\R$. Thus there exists $v$ a solution of the KPP problem \eqref{assKPP}, such that $\ds{u(z)\leq v(z)\leq M}$ for all $z\in\R$. But as $c>\overline{c}$, $v\equiv 0$, which implies that \eqref{elliptproblem} has no positive solution as soon as $c>\overline{c}$ and the Proposition is proved. \fin


\section{Convergence of the Cauchy problem}\label{CauchyPb}
In this section we come back to the parabolic problem \eqref{problem}, that we remind below
\begin{equation*}\begin{cases}
u_t-u_{xx}=f(x-ct,u) & x\in\R, t>0,\\
u(0,x)=u_0(x) & x\in\R,
\end{cases}
\end{equation*}
where $u_0\in H^2(\R)$ is nonnegative, bounded and compactly supported.\\
Letting $z:=x-ct$, $u$ satisfies the following problem
\be\label{parabolicgeneral}\tag{$\tilde{P}$}
\begin{cases}
u_t-u_{zz}-cu_z=f(z,u) &\forall z\in\R, t>0,\\
u(0,z)=u_0(z), &\text{ for all } z\in\R.
\end{cases}
\ee
We know that such a $u$ exists, using sub- and super-solution arguments with assumptions \eqref{f0} and \eqref{ubdd}. As $s\mapsto f(\cdot,s)$ is 
Lipschitz-continuous, applying the maximum principle we have that $u$ is unique. 
Defining $\ds{v(t,z)=u(t,z)e^{\frac{c}{2}z}}$ for all $t>0$, $z\in\R$, then $v$ satisfies the following equation 
$$v_t-v_{zz}+\frac{c^2}{4}v=e^{\frac{c}{2}z}f(z,e^{-\frac{c}{2}z}v).$$
Multiplying the previous equation by $v$ and using \eqref{favbdd} we have
$$\frac{d}{dt}\int_\R\frac{v^2}{2}dz+\int_\R|v_z|^2dz+(\delta+\frac{c^2}{4})\int_\R v^2dz\leq \int_{-R}^R (\frac{f(t,e^{-\frac{c}{2}z}v)}{e^{-\frac{c}{2}z}v}+\delta)v^2dz\leq (||f||_{lip}+\delta)M^2e^{cR}\times 2R,$$
using \eqref{flip} and the fact that $v(\cdot,z)\to 0$ as $z\to\pm\infty$, 
as $v_0$ is bounded and compactly supported and $v$ satisfies $v_t-v_{zz}+(\frac{c^2}{4}+\delta)v\leq0$ when $|z|>R$. 
Proceeding as in \cite{Evans} for example, as $v(0,\cdot)\in H^2(\R)$, we get that $v\in L^2((0,T),H^2(\R))$ and $v_t\in L^2((0,T),L^2(\R))$, 
for all $T>0$. And thus as soon as $u_0\in H^2_c(\R)$, there exists a unique $u\in L^2([0,T[, H_c^2(\R))$, with $u_t\in L^2((0,T),L^2_c(\R))$ for all $T>0$,
solution of \eqref{parabolicgeneral}. Moreover $u(t,z)>0$ for all $t>0$, $z\in\R$. We will now prove Theorem \ref{uniqlimit} on the convergence of solution of \eqref{parabolicgeneral} as $t\to+\infty$. In \cite{Mzero} Matano proves the convergence of solutions of one dimensional semilinear parabolic equations in bounded domains using a geometric argument and the maximum principle and Du and Matano extended this result in \cite{DM} to unbounded domains for homogeneous $f$. Their method relies on classification of solutions for homogeneous problems and uses a reflexion principle which cannot be applied in our case. An alternative proof of this result was first given by Zelenyak in \cite{Zelenyak} using a variational approach. In \cite{HR} Hale and Raugel proved an abstract convergence result in gradient like systems which might apply in the present framework. It roughly states that if the kernel of the linearized equation near any equilibrium has dimension 0 or 1, then the solution of the Cauchy problem converges. We prove such an intermediate step in Lemma \ref{lem:dimL}. We chose to prove directly the convergence of the Cauchy problem in section \ref{Uniqlimitsect} using arguments inspired from Zelenyak's paper \cite{Zelenyak}. But we had to deal with some additional difficulties coming from the fact that our equation is set in $\R$, which induced a lack of compactness and the necessity of finding some controls at infinity. All of this is detailed in section \ref{Uniqlimitsect}. In the next section we start by pointing out the convergence up to a subsequence of the solution $u$ of \eqref{parabolicgeneral}. 


\subsection{Convergence up to a subsequence}\label{convsubseq}
\begin{Prop}\label{convsubsequence} Let $u\in L^2([0,T[, H_c^1(\R))$ for all $T>0$, be the solution of \eqref{parabolicgeneral}. Then there exists a sequence $(t_n)_n$ that goes to infinity as $n\to+\infty$, such that $u(t_n,z)$ converges to a solution of \eqref{elliptproblem} as $n\to+\infty$ locally in $z\in\R$.
\end{Prop}

\noindent\textbf{Proof of Proposition \ref{convsubsequence}:}
As $u(t,\cdot)\in H_c^2(\R)$ and $u_t(t,\cdot)\in L_c^2(\R)$ for all $t>0$, standard arguments show that $t\mapsto E_c[u(t,\cdot)]$ is $C^1$ and
\begin{align*}
\frac{d}{dt}E_c[u(t,\cdot)]&=\int_\R e^{cz}\left\{u_{zt}u_z-f(z,u)u_t\right\}dz\\
			&=\int_\R  (e^{cz}u_z)u_{tz}dz-\int_\R e^{cz}f(z,u)u_tdz\\
			&=-\int_\R (cu_z+u_{zz})e^{cz}u_tdz-\int_\R e^{cz}f(z,u)u_tdz\\
			&=\int_\R (-cu_z-u_{zz}-f(z,u))e^{cz}u_tdz\\
			&=\int_\R -(u_t)^2e^{cz}dz \leq 0.
\end{align*}
We know from Proposition \ref{Eminore} that $E_c[u]$ is bounded from below. It implies that $\ds{E_c[u]\to C} \text{ as } t\to+\infty$, and there exists $(t_n)_n,$ such that $t_n\to+\infty$ and $\ds{\frac{d}{dt}E_c[u](t_n)\to0}$ as $\ds{n\to+\infty}$, i.e $\ds{\lVert u_t(t_n,\cdot)\rVert_{L^2_c(\R)}\to0}$ as $n\to+\infty$, which implies from standard arguments, that up to extraction $u_t(t_n,z)\to0$ as $n\to+\infty$ for almost every $z\in\R$. Using Schauder Theory, we have that $(u(t_n,z))_n$ converges toward $u_\infty$ a stationary solution of \eqref{parabolicgeneral}, i.e a solution of \eqref{elliptproblem}, up to extraction.\fin
\bigskip\\
Now we investigate the uniqueness of the limit $u_\infty$.


\subsection{Uniqueness of the limit}\label{Uniqlimitsect}
We want to prove that, considering compactly supported initial data $u_0$, the solution of our parabolic problem \eqref{parabolicgeneral} admits a unique limit.
Define the $\omega$-limit set:
$$\ds{\Omega(u_0)=\underset{t>0}{\cap}\overline{\left\{u(\tau,\cdot), \quad \tau\geq t\right\}}}.$$
The closure is taken with respect to the topology of $H^2_c(\R)$.

\noindent We first prove the following Lemma,
\begin{Lemma}\label{omegalimitstationnary}
If $w\in \Omega(u_0)$, then $w$ is a solution of the stationary equation $$-w_{zz}-cw_z=f(z,w)\quad \text{in } \R.$$
\end{Lemma}

\noindent\textbf{Proof }: If $w\in\Omega(u_0)$, then there exists a sequence $(t_n)_{n\geq1}$ that converges to $+\infty$ as $n\to+\infty$ such that $u(t_n,z)\to w(z)$ in $H_c^2(\R)$ as $n\to+\infty$. Let $u^n(t,z)=u(t+t_n,z)$ for all $t>0$ and $z\in\R$, then using parabolic estimates, $u^n\to \bar{w}$ as $n\to+\infty$ (up to a subsequence) with $\bar{w}$ solution of \eqref{parabolicgeneral} such that $\bar{w}(0,z)=w(z)$ for all $z\in\R$. Moreover as $E_c[u]$ is decreasing in $t$ and bounded from below $E_c[u^n(t,\cdot)]\to C$ as $n\to+\infty$ and thus $E_c[\bar{w}]=C$ for all $t>0$. We have
\begin{align*}
\frac{d}{dt} E_c[\bar{w}]=0,
\end{align*}
this implies that $\int_\R e^{cz}(\bar{w}_t)^2dz=0$. We thus obtain that $\bar{w}=w$ is a stationary solution of \eqref{parabolicgeneral}, i.e a solution of \eqref{elliptproblem} and we have proved the Lemma. \fin
\bigskip\\

\noindent We will need to prove some Lemmas before starting the proof of the Theorem \ref{uniqlimit}.\\
Take $w \in \Omega (u_0)$. Let 
\begin{equation} \label{defF} \begin{array}{rrcl}
   B:& H^2_c(\R) &\to& L^2_c(\R),\\
&w&\mapsto & w''+cw'+f(z,w).\\ 
  \end{array}\end{equation}
We know that $w$ is a stationary solution of \eqref{elliptproblem}, in other words, $B(w)=0$. 
Define the linear operator:
$$\begin{array}{rrcl}
\L_{w}:= DB (w):& H^2_c(\R) &\to& L^2_c(\R),\\
&h&\mapsto & h''+ch'+f_u'\big(z,w (z)\big)h.\\ 
  \end{array}$$

\begin{Lemma} \label{lem:unifdecay}
Assume that $w\in H^2_c (\R)$ is a non negative, bounded solution of $-w''-cw'\leq -\delta w$ on $\R\backslash (-R,R)$ such that $w(\pm\infty)=0$, then 
$w(z)\leq w(-R) e^{\lambda_+ (z+R)}$ for all $z\leq -R$ and $w(z) \leq w(R) e^{\lambda_- (z-R)}$ for all $z>R$, 
where $\lambda_-<0<\lambda_+$ are the solutions of $\lambda^2+\lambda c = \delta$. 
\end{Lemma}

\noindent \textbf{Proof.}
Define 
$$\phi_-(z):=w (-R)e^{\lambda_+(z+R)}, \quad \forall\: z<-R,$$
$$h(z):=w(z)-\phi(z), \quad \forall\: z<-R.$$
Then $h$ is solution of 
\begin{equation*}
\begin{cases}
-h''-ch'+\delta h\leq 0 &\text{for all } z<-R,\\
h(-\infty)=0,\quad h(-R)\leq 0.
\end{cases}
\end{equation*}
Assume that $h$ achieves a maximum at $z_0\in(-\infty,-R)$. This would imply that $h(z_0)\leq0$ and thus $h\leq0$ in $(-\infty,-R]$. 
Otherwise either $h$ admits a minimum in $(-\infty,-R)$ or is monotone on $(-\infty,-R)$, which also implies that $h\leq0$ in $(-\infty,-R]$ and the first inequality is proved. The inequality on $[R,\infty)$ is proved 
similarly.
\fin

\begin{Lemma} \label{lem:expdecay}
There exists $z_-\in\R$ such that, if $w_1,w_2\in H^2_c (\R)$ are two positive, bounded, solutions of 
$w''+cw'+f(z,w)=0$ over $\R$ with $w_1 (z)=w_2(z)$ for some $z\leq z_-$, then $w_1\equiv w_2$.
\end{Lemma}

\noindent {\bf Proof.}
Let $u:=(w_1-w_2)^2$. This function satisfies
$$\begin{array} {rcl} u''+cu' &=& 2(w_1'-w_2')^2 +2\big( -f(z,w_1) + f(z,w_2)\big)(w_1-w_2)\\
   &\geq& -2f_u'(z,0)u -2\big|- f(z,w_1) + f(z,w_2)-f_u'(z,0)(w_2-w_1) \big||w_1-w_2|. \\
  \end{array}
$$
On the other hand, Lemma \ref{lem:unifdecay} and the $\mathcal{C}^1$ smoothness of $f(z,s)$ with respect to $s$ yields that there exists $z_-$ such that
$$\forall z\leq z_-, \quad |f(z,w_2)-f(z,w_1) - f_u'(z,0) (w_2-w_1)|\leq \frac{\delta}{2} |w_2-w_1|$$
where $\delta$ is the constant defined by \eqref{favbdd}. We thus get
$$\forall z\leq z_-, \quad u''+cu' \geq -2f_u'(z,0)u -\delta u \geq \delta u$$
decreasing $z_-$ once more if necessary. 

It now follows from this inequation that $u$ cannot reach any local maximum over $(-\infty,z_-)$. As $u(-\infty)=0$ and $u\geq 0$, it implies that $u$ is nondecreasing. 
Lastly, if $w_1 (z) = w_2 (z)$ for some $z\leq z_-$, then $u(z)=0$ and thus $u\equiv 0$, meaning that $w_1\equiv w_2$. 
\fin

\begin{Lemma} \label{lem:dimL}
 $$ \mathrm{ dim Ker} \L_{w} \in \{0,1\}.$$
\end{Lemma}

\noindent {\bf Proof.} The Cauchy theorem yields that
$$\mathrm{ Ker} \L_{w} = \{ h\in H^2_c (\R), \quad h''+ch'+f_u'\big( z, w (z)\big)h=0\}$$
has at most dimension $2$. If it has dimension $2$, then it would mean that for all $z_0\in\R$ and for all couple $(h_0,h_1)$, the solution of 
$$h''+ch'+f_u'\big( z,w(z)\big)h=0, \quad h(z_0)=h_0, \quad h'(z_0)=h_1$$
belongs to $H^2_c(\R)$. In particular $h(+\infty)=0$.

But now the same arguments as in the proof of Lemma \ref{lem:expdecay} yields that $h^2$ is nonincreasing over $(z_+,+\infty)$ and thus one reaches a contradiction by
taking $z_0>z_+$ and $(h_0,h_1)$ such that $h_0 h_1>0$. 
\fin

\begin{Lemma} \label{lem:Luf}
Assume that $\mathrm{ dim Ker} \L_{w} =1$. Then there exists a constant $C=C(w)$ such that for all $g\in L^2_c (\R)$,
if $u\in H^2_c (\R)$ satisfies $\L_w u = g$ in $\R$ and $\int_\R e^{cz} u(z) v(z)dz=0$ for all $v\in \mathrm{ Ker} \L_{w}$, then 
$$ \|u\|_{H^2_c (\R)} \leq C \|g\|_{L^2_c (\R)}.$$
Moreover, if $W$ is a family of solutions $w\in H^1_c(\R)$ of \eqref{elliptproblem} such that $\mathrm{ dim Ker} \L_{w} =1$ for all $w\in W$ and $\sup_{w\in W} \|w\|_{H^1_c(\R)}<\infty$, then 
the constant $C$ can be chosen to be the same for all $w\in W$. 
\end{Lemma}

\noindent {\bf Proof.}
Clearly the operator 
$$\begin{array}{rrcl}
T:& (\mathrm{ Ker}\L_w)^{\perp} &\to& \mathrm{Im} \L_w\\
&h&\mapsto & \L_w h \\ 
  \end{array}$$
where $(\mathrm{ Ker}\L_w)^{\perp}$ is with respect to the scalar product of $L_c^2$. $T$ is invertible and continuous. Hence the bounded inverse theorem yields that its inverse is continuous. Taking $C$ its continuity constant, this means that 
for all $g\in L^2_c (\R)$ such that there exists $u\in (\mathrm{Ker} \L_w)^{\perp}$ satisfying $\L_w u = g$, one has 
$\|u\|_{H^2_c(\R)} \leq C \|g\|_{L^2_c (\R)}$ and the result follows. 
\bigskip\\
Next, we first prove that there exists $C>0$ such that if $W$ is a family of solutions $w\in H_c^1(\R)$ of \eqref{elliptproblem} such that $\mathrm{Ker} \L_w \neq \{0\}$ for all $w\in W$ and $\sup_{w\in W} \|w\|_{H^1_c(\R)}<\infty$, then 
$$ \|u'\|_{L^2_c (\R)} \leq C \|g\|_{L^2_c (\R)}.$$
Assume that this is not true, there would exist a sequence $(w_n)_n$ of solutions of \eqref{elliptproblem}, bounded in $H^1_c(\R)$, such that 
$\mathrm{Ker} \L_{w_n} \neq \{0\}$ for all $n$ and the associated constants $C_n=C(w_n)$ converge to $+\infty$ as $n\to +\infty$. 
In other words, there exist $v_n\in \mathrm{Ker} \L_{w_n}$ for all $n$ and two sequences $(u_n)_n$ in $H^2_c(\R)$ and $(g_n)_n$ in $L^2_c(\R)$ 
such that $\L_{w_n}u_n = g_n$ in $\R$, $\int_\R e^{cz} u_n(z)v_n(z)dz=0$, $\|u_n'\|_{L^2_c(\R)}=1$ for all $n$ and $\lim_{n\to +\infty} \|g_n\|_{L^2_c(\R)}=0$. 
Up to multiplication, we can assume that $\|v_n'\|_{L^2_c(\R)}=1$. \\
As $(w_n)_n$ is bounded in $H^1_c(\R)$, we can assume, up to extraction, that it converges locally uniformly to some function $w_\infty \in H^1_c(\R)$. 
Similarly, the Poincar\'e inequality stated in Lemma \ref{lemmapoincareH1c} yields that $(u_n)_n$ and $(v_n)_n$ are indeed bounded in $H^1_c(\R)$ and 
we can thus define their weak limits $u_\infty$ and $v_\infty$ in $H^1_c(\R)$. 
As $u_n'' = -cu_n' -f_u'\big(z, w_n(z)\big) u_n +g_n$, multiplying by $u_ne^{cz}$ and integrating over $\R$, as $\|u_n'\|_{L^2{_c}(\R)}=1,$ we get
$$1-\int_\R e^{cz}f'_u(z,w_n)u_n^2dz=-\int_\R e^{cz}u_ng_ndz.$$
As $u_n$ converge weakly in $L_c^2$ and $g_n\to0$ in $L_c^2$ as $n\to+\infty$, the right-hand side converges to $0$ as $n\to+\infty$. 
Assuming $u_n\rightharpoonup 0$ in $L^2_c$ yields a contradiction. Indeed, using Lemma \ref{lem:unifdecay} for all $n$, for all $\epsilon>0$ there exists $r>0$ such that $w_n(z)<\epsilon$ for all $|z|>r$. As $f(z,\cdot)$ is $C^1$, for $\epsilon$ small enough, $f'_u(z,w_n)<0$ for all $|z|>r$. And we obtain
$$1-\int_{-r}^re^{cz}f'_u(z,w_n)u_n^2dz\leq -\int_\R e^{cz}u_ng_ndz,$$
which yields a contradiction when we let $n\to+\infty$, as $u_n\to0$, strongly in $L^2_{c}([-r,r])$. This implies that $u_\infty\not\equiv 0$.\\
Using the same arguments with $v_n$, as $\|v_n'\|_{L^2_c(\R)}=1$ for all $n$, we have that
$$1-\int_\R e^{cz}f'_u(z,w_n)v_n^2dz=0.$$
On the other hand, it follows from Lemma \ref{lem:unifdecay} that one can apply the dominated convergence theorem using the bounds 
$v_n (z) \leq M$ for all $z<R$ and $v_n (z) \leq Me^{\lambda_- (z-R)}$ for all $z>R$, since $\ds{c<-2\lambda_-}$. We thus obtain
$$1-\int_\R e^{cz}f'_u(z,w_\infty)v_\infty^2 dz=0.$$
Moreover, classical elliptic regularity estimates yield that $v_\infty$ satisfies $\L_{w_\infty}v_\infty=0$ in $\R$. Integrating by parts, we get
$$\int_\R e^{cz} \Big\{ (v_\infty')^2 - f_u'(z,w_\infty) v_\infty^2\Big\} dz = 0.$$
We thus conclude that $\|v_\infty'\|_{L^2_c (\R)}=1$. As $L^2_c (\R)$ is an Hilbert space, this indeed implies that 
$(v_n')_n$ converges strongly to $v_\infty'$ in $L^2_c (\R)$ as $n\to +\infty$. 
Using the Poincar\'e type inequality given in Lemma \ref{lemmapoincareH1c} we have that $v_n\to v_\infty$ in $L^2_c$ as $n\to+\infty$. 
This implies that $\int_\R e^{cz} u_\infty (z)v_\infty (z)dz=0$, $\L_{w_\infty} u_\infty = 0$ and $\L_{w_\infty} v_\infty=0$ 
over $\R$. Hence, $\mathrm{ dim Ker} \L_{w_\infty} = 2$, which contradicts Lemma \ref{lem:dimL}. 

Thus there exists $C>0$ such that if $W$ is a family of solutions $w\in H_c^1(\R)$ of \eqref{elliptproblem} such that $\mathrm{Ker} \L_w \neq \{0\}$ for all $w\in W$ and $\sup_{w\in W} \|w\|_{H^1_c(\R)}<\infty$, then 
$$ \|u'\|_{L^2_c (\R)} \leq C \|g\|_{L^2_c (\R)}.$$
Now to prove the last assertion of the Lemma we just apply Lemma \ref{lemmapoincareH1c} to get a bound on $||u||_{L_c^2(\R)}$ and the bound on $||u''||_{L_c^2(\R)}$ follows from the equation.
\fin

\begin{Lemma}  \label{lem:lem4}
Assume that for some $T>0$, there exist two constants $K,C>0$ such that for all $t\in [0,T]$,
$$\int_t^\infty\int_\R e^{cz}  u_t^2 (s,z)dsdz \leq K e^{-Ct}.$$
Then for all $0\leq t\leq \tau\leq T$, one has:
$$|| u(t,\cdot)-u(\tau,\cdot)||_{L^2_c (\R)} \leq  \frac{\sqrt{K}}{1- e^{-C/2}} e^{-Ct/2}.$$
\end{Lemma}

\noindent {\bf Proof:} This Lemma is similar to Lemma 4 in Zelenyak paper \cite[Lemma 4]{Zelenyak}. As our solutions are defined on the full line $\R$ instead of a segment, we obtain a control in $L^2$ instead of $L^1$.\\
Assume first that $|t-\tau|\leq 1$. Then 
$$\begin{array}{rcl}
  || u(t,\cdot)-u(\tau,\cdot)||^2_{L^2_c (\R)}&=& \int_\R e^{cz} \Big|\int_t^\tau u_t (s,z) ds\Big|^2 dz\\
&&\\
&\leq &  \int_\R \int_t^\tau  (\tau-t) e^{cz} u_t^2 (s,z) ds dz\\
&&\\
&\leq & \int_\R \int_t^\infty  e^{cz} u_t^2 (s,z) ds dz\\
&&\\
&\leq & K e^{-Ct}.\\
 \end{array}$$

Next, if $|t-\tau|> 1$, let $N=[\tau-t]$ be the integer part of $\tau-t$. We compute:
$$\begin{array}{rcl}
  || u(t,\cdot)-u(\tau,\cdot)||_{L^2_c (\R)}&\leq &  \sum_{n=0}^{N-1}|| u(t+n,\cdot)-u(t+n+1,\cdot)||_{L^2_c (\R)}+  || u(t+N,\cdot)-u(\tau,\cdot)||_{L^2_c (\R)}  \\
&&\\
&\leq & \sum_{n=0}^{N-1} \sqrt{K} e^{-C (t+n)/2}+  \sqrt{K} e^{-C(t+N)/2}\\
&&\\
&\leq &  \ds\frac{\sqrt{K}}{1- e^{-C/2}}e^{-Ct/2},  \\
 \end{array}$$
which ends the proof. \fin
\bigskip\\
\noindent \textbf{Proof of Theorem \ref{uniqlimit}:} 
Assume that $\Omega(u_0)$ is not an isolated point. Using Lemma \ref{omegalimitstationnary} and Lemma \ref{lem:expdecay} we can choose $R$ large enough such that $\Omega(u_0)$ 
is parametrised by the value of the function at $-R$, i.e $\ds{\Omega(u_0)=\{w(\alpha, \cdot), \quad w(\alpha,-R)=\alpha 
\text{ and } w \text{ is a stationary solution}\}}$. As $u$ is bounded, the quantities
$\ds{0\leq \alpha_1=\underset{t\to+\infty}{\liminf}u(t,-R)<\alpha_2=\underset{t\to+\infty}{\limsup}u(t,-R)}$ are well-defined and 
classical connectedness and compactness arguments yield that for all $\alpha\in[\alpha_1,\alpha_2]$, as $t\mapsto u(t,-R)$ is continuous, by Schauder parabolic regularity estimates, there exists $(t_n)_n$ such that $u(t_n,-R)\to\alpha$ as $n\to+\infty$. Hence $\Omega (u_0)$ is the curve $\{ w(\alpha,\cdot), \alpha \in [\alpha_1,\alpha_2]\}$. 

For each $w(\alpha,\cdot)\in\Omega(u_0)$, we have that $\alpha\mapsto w(\alpha,\cdot)$ is continuous with respect to the 
$L^\infty (\R)$ norm, using the uniqueness property proved in Lemma \ref{lem:expdecay} and the elliptic regularity. Define 
$$v_\epsilon:=\frac{w(\alpha+\epsilon,\cdot)-w(\alpha,\cdot)}{\epsilon},$$
it satisfies 
$$-v''_\epsilon-cv_\epsilon'-(f'_u(z,w(\alpha,\cdot))+R_\epsilon)v_\epsilon=0,$$
with $\|R_\epsilon\|_{L^\infty (\R)}\to0$ as $\epsilon\to0$ since $\alpha\mapsto w(\alpha,\cdot)$ is continuous and $s\mapsto f(\cdot,s)$ is $C^1$. 
Then applying Lemma \ref{lem:unifdecay} to the nonnegative and nonpositive parts of $v_\epsilon$, we get that 
$$|v_\epsilon(z)|\leq |v_\epsilon(-R)|e^{\lambda_+(z+R)}=e^{\lambda_+(z+R)},\quad \forall\:z\in(-\infty,-R),$$
where we have used that 
$$v_\epsilon (-R)=\frac{w(\alpha+\epsilon,-R)-w(\alpha,-R)}{\epsilon}=\frac{\alpha+\epsilon-\alpha}{\epsilon}.$$
This implies for all $h>0$:
$$\frac{e^{-\lambda_+h}-1}{h}\leq\lvert\frac{v_\epsilon(-R-h)-v_\epsilon(\alpha,-R)}{-h}\rvert\leq \frac{-e^{-\lambda_+h}+1}{h}.$$
Letting $h\to0$ we get that $v_\epsilon'(-R)$ is bounded uniformly with respect to $\epsilon$. 
The continuity with respect to the initial condition and parameters for ordinary differential equations yields that 
$(v_\epsilon)_\epsilon$ is bounded in $C^2_{loc}$ and using Lemma \ref{lem:unifdecay} again
$$|v_\epsilon(z)|\leq |v_\epsilon(R)|e^{-\lambda_-(z-R)}, \quad \forall\:z\in(R,+\infty).$$
Up to a subsequence $v_\epsilon\to v$ in $C^2_{loc}$ by elliptic regularity, and $v\in L_c^2(\R)$ using the uniform exponential convergence.
We have thus proved that
$\ds{v=\frac{\partial w}{\partial \alpha}}$ is well-defined and is a solution of 
$$v(-R)=1 \quad \hbox{and} \quad \L_{w} v=v''+cv'+f_u'\big(z,w(\alpha,\cdot)\big)v=0 \quad \hbox{ over } \R.$$
Lastly, multiplying the equation by $ve^{cz}$, integrating over $\R$ we have that
$$\int_\R (v_z)^2e^{cz}dz=\int_\R f'_u(z,w(\alpha,\cdot))v^2e^{cz}dz,$$
integrating by part and using the exponential convergence of $v$. As $f(\cdot, s)$ is non positive outside $(-R,R)$, we get
$$\int_\R(v_z)^2e^{cz}dz \leq C\int_{-R}^Rv^2e^{cz}dz.$$
Hence $v\in H^1_c(\R)$.\\
\noindent We have $v\not\equiv 0$ in $\R$. Now we define for fixed $t>0$,
$$\alpha(t):=\arg\inf\left\{||u(t,\cdot)-w(\alpha,\cdot)||_{L^2_c(\R)}, \alpha\in[\alpha_1,\alpha_2]\right\}.$$
For each $t>0$, if the $\inf$ is attained at an interior point $\alpha(t)\in (\alpha_1,\alpha_2)$, then 
$\ds{\frac{\partial}{\partial\alpha}||u(t,\cdot)-w(\alpha,\cdot)||^2_{L^2_c(\R)}\Big|_{\alpha=\alpha(t)}}=0$, and thus
$$\int_\R e^{cz}\big(u(t,z)-w(\alpha,z)\big)\frac{\partial w}{\partial\alpha}\big|_{\alpha=\alpha(t)}dz=0.$$
We thus have for all $t>0$ such that $\alpha (t) \in (\alpha_1,\alpha_2)$: 
$$\L_{w (\alpha (t),\cdot)} v=0, \quad \int_\R e^{cz}(u-w)v |_{\alpha=\alpha(t)} dz=0\quad \hbox{ and } \quad \L_{w (\alpha (t),\cdot)}  (u-w)=g,$$ 
with 
$$g(t,z):=u_t(t,z)+b(t,z)(u(t,z)-w\big(\alpha (t),z)\big),$$
$$b(t,z):=f_u'\big(z,w(\alpha (t),z)\big)-\frac{f\big(z,u(t,z)\big)-f\big(z,w(\alpha (t),z)\big)}{u(t,z)-w( \alpha (t),z)}.$$
Lemma \ref{lem:Luf} thus applies and gives
$$||u(t,\cdot)-w(\alpha(t),\cdot)||_{H^2_c(\R)}\leq C||u_t(t,\cdot)||_{L^2_c(\R)}+ C||b(t,\cdot)||_{L^2_c(\R)}  ||u(t,\cdot)-w(\alpha(t),\cdot)||_{L^2_c(\R)} ,$$
for all $t>0$ such that $\alpha(t)\in (\alpha_1,\alpha_2)$.
But as $f=f(z,u)$ is of class $\mathcal{C}^1$ with respect to $u$ uniformly in $z$ and as $\lim_{t\to +\infty}||u(t,\cdot)-w(\alpha (t),\cdot)||_{L^2_c(\R)}=0$, one has 
$||b(t,\cdot)||_{L^2_c(\R)}\to 0$ as $t\to+\infty$ and it thus follows that, even if it means increasing $C$, for all admissible $t>0$, one has 
$$||u(t,\cdot)-w(\alpha(t),\cdot)||_{H^2_c(\R)}\leq C||u_t(t,\cdot)||_{L^2_c(\R)}.$$
and $C$ is bounded independently of $\alpha (t) \in (\alpha_1,\alpha_2)$.

Now ending the proof as in Zelenyak \cite{Zelenyak}, we have that for all $t>0$ and any $w\in\Omega (u_0)$, the solution $u$ of our parabolic problem satisfies
\begin{align*}
E_c[u(t,\cdot)]-E_c[w]&=\frac{1}{2}\int_\R e^{cz}\big(u_z^2(t,z)-w^2_z (z)\big) dz-\int_\R e^{cz}\big(F(z,u(t,z))-F(z,w(z))\big)dz\\
&= \frac{1}{2}\int_\R e^{cz}(u_z-w_z)^2 dz +\int_\R e^{cz}(u_z-w_z)w_z dz \\
&\phantom{==}-\int_\R e^{cz} f\big( z, w(z)\big) (u(t,z)-w(z))dz +\int_\R e^{cz} \sigma(t,z)(u(t,z)-w(z))^2 dz,\\
\end{align*}
where $\sigma=\sigma(t,z)$ is a bounded and measurable function since $f=f(z,u)$ is of class $\mathcal{C}^1$ with respect to $u$, uniformly in $z$. 
As $w$ is a stationary solution of \eqref{elliptproblem}, integrating by parts, we get
$$\begin{array}{rcl}
E_c[u(t,\cdot)]-E_c[w]&=& \frac{1}{2}\int_\R e^{cz}(u_z-w_z)^2 dz  +\int_\R e^{cz} \sigma(t,z)(u(t,z)-w(z))^2 dz\\
&&\\
&\leq &\sup \{\frac{1}{2}, \|\sigma\|_{L^\infty (\R)}\} ||u(t,\cdot)-w||^2_{H_c^1(\R)}.\\
\end{array}$$

Next, we have shown in the proof of Proposition \ref{convsubseq} that the energy is decreasing and bounded from below, 
so $E_c^\infty := \lim_{t\to+\infty} E_c [u(t,\cdot)]$ is well-defined. 
For all $t>0$ such that $\alpha(t)\in (\alpha_1,\alpha_2)$, gathering the previous inequalities, one gets
\be\label{dtenergy}
\frac{d}{dt}(E_c(u(t,\cdot)-E_c^\infty)=-||u_t(t,\cdot)||^2_{L^2_c(\R)}\leq -C^{-1}||u(t,\cdot)-w(\alpha (t),\cdot)||^2_{H^2_c(\R)}\leq -K(E_c[u(t,\cdot)]-E_c^\infty),
\ee
where $K$ is an explicit constant and $E_c^\infty = \lim_{t\to+\infty} E_c [u(t,\cdot)]$ is equal to $E_c [w]$ for all $w\in \Omega [u_0]$. 

Now let $\alpha_0\in (\alpha_1,\alpha_2)$ and take a sequence $(t_n)_n$ such that $\lim_{n\to +\infty} t_n=+\infty$ and $\lim_{n\to +\infty} u(t_n,z)= w(\alpha_0,z)$. 
There exists $\eta >0$ such that 
$$||w(\alpha_0,\cdot)-w(\alpha_1,\cdot)||_{L_c^2(\R)}>\eta \quad \hbox{ and } \quad ||w(\alpha_0,\cdot)-w(\alpha_2,\cdot)||_{L_c^2(\R)}>\eta.$$
Choose $N$ large enough such that $||u(t_N,\cdot)-w(\alpha_0,\cdot)||_{L_c^2(\R)}\leq \frac{\eta}{8}$ and for all $t\geq t_N$
$$\sqrt{E_c[u(t,\cdot)]-E_c^\infty}\leq (1-e^{-C/2})\frac{\eta}{8}.$$
We set 
$$\bar{t}=\inf\big\{t\geq t_N,\quad ||u(t,\cdot)-w(\alpha_0,\cdot)||_{L_c^2(\R)}\geq \min \{||u(t,\cdot)-w(\alpha_1,\cdot)||_{L_c^2(\R)}, ||u(t,\cdot)-w(\alpha_2,\cdot)||_{L_c^2(\R)}\}\big\} .$$
Clearly $\alpha (t)\neq \alpha_1$ and $\alpha (t)\neq \alpha_2$, that is, $\alpha (t)$ is an interior point, for all $t\in [t_N,\bar{t})$. Hence, 
inequality \eqref{dtenergy} holds for all $t\in[t_N,\bar{t})$, i.e
$$E_c[u(t,\cdot)]-E_c^\infty\leq \left(E_c[u(t_N,\cdot)]-E_c^\infty\right)e^{-C(t-t_N)}.$$
By Lemma \ref{lem:lem4}, one has for all $t_N \leq t\leq \tau\leq \bar{t}$:
\be \label{eq:cvexp} ||u(t,z)-u(\tau,z)||_{L^2_c (\R)}\leq \ds\frac{\sqrt{E_c[u(t_N,\cdot)]-E_c^\infty}}{1-e^{-C/2}}e^{-C(t-t_N)/2}\leq \frac{\eta}{8}e^{-C(t-t_N)/2}.\ee
If $\bar{t}$ is finite then from the previous inequality we obtain that 
\be\label{ineq1}
||u(\bar{t},\cdot)-w(\alpha_0,\cdot)||_{L_c^2(\R)}\leq ||u(\bar{t},\cdot)-u(t_N,\cdot)||_{L_c^2(\R)}+||u(t_N,\cdot)-w(\alpha_0,\cdot)||_{L_c^2(\R)}\leq \frac{\eta}{4}
\ee
and, for $k=1$ and $k=2$:
\be\label{ineq2}
||u(\bar{t},\cdot)-w(\alpha_k,\cdot)||_{L_c^2(\R)}\geq ||w(\alpha_k,\cdot)-w(\alpha_0,\cdot)||_{L_c^2(\R)}-||u(\bar{t},\cdot)-w(\alpha_0,\cdot)||_{L_c^2(\R)}\geq \eta-\frac{\eta}{4}=\frac{3}{4}\eta.
\ee
Comparing \eqref{ineq1} and \eqref{ineq2} we conclude that $\inf ||u(\bar{t},\cdot)-w(\alpha,\cdot)||_{L_c^2(\R)}$ cannot be attained for $\alpha=\alpha_k$,($k=1,2$), 
and thus $\bar{t}=\infty$. We thus conclude that (\ref{eq:cvexp}) holds for all $\tau\geq t\geq t_N$
which proves that $t\mapsto u(t,\cdot)$ converges strongly in $L^2_c$. This contradicts the assumption that $\Omega (u_0)$ is not an isolated point and concludes the proof.
\fin


\section{On the stability of the trivial steady state 0}\label{stab0}
In this section we discuss the different behaviours of the solution of \eqref{problembis} depending on the stability of 0 and the initial condition $u_0$. We first define what we mean by stability of the trivial steady state 0. 
\bigskip\\
Let $\L$ be the linearized operator around 0: 
$$-\L u:=-u''-cu'-f_s(z,0)u,$$
defined for all $u\in H^1(\R)$. It is easy to check (using Lemma \ref{lem:unifdecay}) that the operator $\L$ admits a principal eigenfunction in $H_c^1(\R)$, that is there exist $\ds{(\lambda_c,\phi)}$ such that 
\be\label{L}\begin{cases}
-\L \phi=\lambda_c\phi &\text{in } \R,\\
\phi>0 &\text{in } \R,\\ \phi\in H_c^1(\R).
\end{cases}
\ee
This eigenvalue $\lambda_c$ is also characterized as the generalized eigenvalue of $\L$:
\be\label{geneigen}
\lambda_c(-\L,\R):=\sup\left\{ \lambda\in\R, \hspace{0.1cm} \exists \phi\in W^{2,1}_\text{loc}(\R), \hspace{0.1cm} \phi>0, (\L+\lambda)\phi\leq 0 \text{ a.e in } \R\right\}.
\ee
One can look at \cite{BR1} and references therein for more details about generalized eigenvalue. We know from \cite[Proposition 1 - section 2]{BR1} that, if we denote by $\lambda(r)$ the principal eigenvalue of our problem on $B_r$ with Dirichlet boundary condition, then $\lambda(r)\to \lambda_c$ as $r\to+\infty$ and there exists $\phi_c\in W^{2,p}_\text{loc}(\R)$, $1\leq p<+\infty$, the principal eigenfunction solution of \eqref{L}.
\bigskip\\
Letting $\ds{v(z)=u(z)e^{\frac{cz}{2}}}$, then 
$$-\L u=0 \iff -\tilde{\L}v=-v''+\frac{c^2}{4}v-f_s(z,0)v=0,$$
where $\tilde{\L}$ is self adjoint. From \cite{BR1, BRprep} 
\be\label{eigenvalue}
\lambda_c(-\L,\R)=\lambda_c(-\tilde{\L},\R)=\underset{\phi\in H^1(\R), \phi\not\equiv0}{\inf}\frac{\int_\R \phi'(z)^2+(\frac{c^2}{4}-f_s(z,0))\phi(z)^2dz}{\int_\R \phi(z)^2dz}.
\ee
If we define $\lambda_0$ as the generalized eigenvalue corresponding to $c=0$, i.e when the medium does not move with time, then we have that 
$$\lambda_c=\lambda_0+\frac{c^2}{4}.$$
We will say that $0$ is linearly stable (respectively unstable) if $\lambda_c>0$ (respectively $\lambda_c<0$). Let us notice that if $0$ is stable in the steady frame, i.e $\lambda_0>0$, then 0 is necessarily stable in the moving frame.


\subsection{Convergence to a travelling wave solution when $0$ is linearly unstable}
In this section we want to prove that when 0 is linearly unstable, i.e $\lambda_0<0$ and $c<2\sqrt{-\lambda_0}$, for $u_0\not\equiv0$ nonnegative initial condition, the solution $u$ of \eqref{parabolicgeneral} converges to a non trivial travelling wave solution as time goes to infinity. 
\begin{Prop}\label{speedinterval}
Let us assume that $\lambda_0<0$ and that $f$ satisfies \eqref{f0}-\eqref{favbdd}, then for all $\ds{c<2\sqrt{-\lambda_0}}$,
$$\underset{u\in H^1_{c}(\R)}{\inf}E_{c}[u]<0,$$
i.e there exists a non trivial solution of \eqref{elliptproblem}
\end{Prop}
If $f$ is of KPP-type, that is $u\mapsto\frac{f(z,u)}{u}$ is decreasing, this result is optimal, i.e $\underset{u\in H^1_{c}(\R)}{\inf}E_{c}[u]\geq0$ and travelling wave solutions do not exist if $c\geq2\sqrt{-\lambda_0}$ (cf \cite{BDNZ}). This is not true for general $f$ (see the next section).
\vspace{0.2cm}\\
\textbf{Proof of Proposition \ref{speedinterval}}: 
 Take $\lambda$ such that $\lambda_0<\lambda< -c^2/4$. It follows from (\ref{eigenvalue}) that there exists $\phi_0 \in H^1 (\R)$ such that 
$$\int_\R \big( \phi_0'(z)^2 - f_s (z,0)\phi_0^2(z)\big) dz \leq \lambda \int_\R \phi_0^2(z)dz.$$
Let $$\phi_n(z)=\frac{\phi_0(z)e^{-\frac{c}{2}z}}{n} \hspace{0.2cm} \forall z\in\R.$$
Then we have the following computation:
\begin{align*}
E_c[\phi_n]&=\int_\R\left\{\frac{|(\phi_0(z)e^{-\frac{c}{2}z})_z|^2}{2n^2}-F\left(z,\frac{\phi_0(z)e^{-\frac{c}{2}z}}{n}\right)\right\}e^{cz}dz\\
		&=\int_\R \frac{(\phi_0'(z))^2}{2n^2}+\frac{c^2}{4}\frac{(\phi_0(z))^2}{2n^2}\\
		&\hspace{0.8cm}-\left(F(z,0)+F_s(z,0)\frac{\phi_0(z)e^{-\frac{c}{2}z}}{n}+F_{ss}(z,0)\frac{(\phi_0(z)e^{-\frac{c}{2}z})^2}{2n^2}+o\Big(\frac{(\phi_0(z)e^{-\frac{c}{2}z})^2}{n^2}\Big)\right)e^{cz}dz\\
		&=\int_\R \frac{(\phi_0'(z))^2}{2n^2}+\frac{c^2}{4}\frac{(\phi_0(z))^2}{2n^2}\\
		&\hspace{0.8cm}-\left(f(z,0)\frac{\phi_0(z)e^{-\frac{c}{2}z}}{n}+f_{s}(z,0)\frac{(\phi_0(z)e^{-\frac{c}{2}z})^2}{2n^2}+o\Big(\frac{(\phi_0(z)e^{-\frac{c}{2}z})^2}{n^2}\Big)\right)e^{cz}dz\\
		&\leq \int_\R (\lambda+\frac{c^2}{4})\frac{(\phi_0(z))^2}{2n^2}dz+o(\frac{1}{n^2}).
\end{align*}
This implies that 
$$\underset{u\in H_c^1(\R)}{\min} E_c[u]\leq(\lambda+\frac{c^2}{4}) \int_\R \frac{(\phi_0(z))^2}{2n^2}dz+o(\frac{1}{n^2})<0,$$
for $n$ large enough. The Proposition is proved. \fin
\bigskip\\
And we have the following Proposition to characterize the behaviour of $u$ as time goes to infinity.

\begin{Prop}\label{0unstable}
If $\lambda_0<0$, for all $c<2\sqrt{-\lambda_0}$, the solution $u$ of \eqref{parabolicgeneral} converges to a non trivial solution of \eqref{elliptproblem} as $t\to+\infty$.
\end{Prop}
\textbf{Proof of Proposition \ref{0unstable}}: We will use the same argument as in \cite[section 2.4]{BR1}. We know that $\lambda(R)\to\lambda_c$ as $R\to+\infty$, and $\lambda_c<0$ thus for $R$ large enough $\lambda(R)<0$ and let $\phi_R>0$ be the principal eigenfunction. Define 
\be\label{Usubsol}
\underline{U}=\begin{cases} \kappa \phi_R &\text{in } (-R,R),\\
0 &\text{otherwise},
\end{cases}
\ee
Then for $\kappa$ small $\underline{U}$ is a subsolution of \eqref{parabolicgeneral} and $\underline{U}\leq u(\tau, \cdot)$ in $\R$ for some $\tau>0$ small, 
$\overline{U}\equiv M\geq u_0$ in $\R$ and is a super solution. Then the solution u of \eqref{parabolicgeneral} is greater than $\underline{U}$ for all $t>0$ and $x\in\R$. Moreover using Theorem \ref{uniqlimit} we know that $u$ converges to $u_\infty\geq \underline{U}$ as $t\to+\infty$. And thus $u$ converges to a non trivial travelling wave solution as $t\to+\infty$.\fin

\begin{Rk}
Notice that Proposition \ref{0unstable} contains Proposition \ref{speedinterval} but the proof of Proposition \ref{speedinterval} 
is interesting as it exhibits the link between the energy and the principal eigenvalue.
\end{Rk}

\subsection{Existence of a travelling wave with positive energy when $0$ is linearly stable}  \label{sec:bistable}
In this section we use the same notations than in the previous one and assume now that 
\be\label{0stableeq}
\lambda_c>0 \quad \hbox{and} \quad \inf_{u\in H^1_c(\R)} E_c [u]<0.
\ee
In this framework, we show that the Mountain Pass Theorem applies and gives the existence of a travelling wave solution with positive energy. 
This provides a class of examples for which travelling wave solutions are not unique. 
Moreover, we construct in section \ref{nonuniqprofsection} an example where uniqueness does not hold even in the class of stable travelling wave solutions with negative energy.\\
We also exhibit at the end of this section the dependence of the asymptotic limit on the initial condition 
and we show in Proposition \ref{0stable} that, depending on the initial condition, we can converge either to 0 or to a travelling wave solution.

\begin{Prop}\label{enerpos}
Assume that $\lambda_c>0$, if $\ds{\underset{u\in H_c^1(\R)}{\min} E_c[u]<0}$ then there exists at least two non trivial travelling wave solution of \eqref{elliptproblem} and one of them has a positive energy.
\end{Prop}

An easy application of this proposition is the following Corollary. 

\begin{Cor}\label{critpoint}
Let 
\begin{equation*}
f(z,u)=\begin{cases} f_0(u) & \text{if } |z|<R,\\
				-\delta u & \text{otherwise},
	\end{cases}
\end{equation*}
where $R,\delta>0$, $f_0$ is a bistable function., i.e 
\begin{gather*}\label{bistable}
\text{There exists } \theta\in (0,1) \text{ such that } f_0(0)=f_0(\theta)=f_0(1)=0,\text{ and } f_0'(0)<0,\hspace{0.3cm} f_0'(1)<0,\\
f_0(s)<0 \text{ for all } s\in(0,\theta)\textrm{, }f_0(s)>0 \text{ for all } s\in (\theta,1),
\end{gather*}
with positive mass:
\begin{equation}\label{intf}
\int_0^1f_0(\tau)d\tau>0.
\end{equation}
Then for $R$ sufficiently large, there exists $\underline{u},\: \overline{u}\in H_c^1 (\R)$ solution of \eqref{elliptproblem} such that $E_c[\underline{u}]<0$ and $E_c[\overline{u}]>0$.
\end{Cor}

Let us highlight this corollary which is totally different from what is known when $f$ satisfies the KPP property. 
Indeed in the present framework 0 is linearly stable, nevertheless we still have the existence of travelling wave solutions. 
\bigskip\\
\noindent {\bf Proof of Corollary \ref{critpoint}.} As $f_s(z,0) = f_0'(0)<0$ if $|z|<R$, $-\delta<0$ otherwise, one has $\lambda_0>0$ and thus $\lambda_c = \lambda_0+c^2/4>0$. 

Moreover, as $f_0$ has a positive mass, taking 
\be
u_{min}(z)=\begin{cases}1 &\text{for all } |z|<R,\\
					0 &\text{for all } |z|>R+1,
				\end{cases}
\ee
such that $u_{min}\in H^1_c(\R)$, one can check that for $R$ large enough $E_c[u_{min}]<0$ and $\ds{\lVert u_{min}\rVert_{H_c^1}>r}$. Proposition \ref{enerpos} applies and gives the conclusion. 
\fin

\bigskip

To prove Proposition \ref{enerpos} we start with the following Lemma.
\begin{Lemma}\label{pointhaut}
For all $r>0$ small enough, one has $\inf_{\lVert u\rVert_{H_c^1 (\R)}=r} E_c [u]>0$. 
\end{Lemma}

\noindent {\bf Proof of Lemma \ref{pointhaut}}: To prove this Lemma, we just need to prove that 0 achieves a strict local minimum, i.e $dE_c[0]\equiv 0$ and $d^2E_c[0]>0$ in the sense that for all $w\in H_c^1(\R)$, $w\not\equiv0$, $d^2E_c[0](w,w)>0$, with 
$$d^2 E_c[0](w,w)=\int_\R e^{cz}\left\{w_z^2-f_s(z,0)w^2\right\}dz.$$
Using the equalities in \eqref{eigenvalue} with $\phi (z) = e^{c z/2} w(z)$, we get,
$$d^2E_c[0](w,w)\geq \lambda_c \lVert w\rVert^2_{H_c^1(\R)},$$
for all $w\in H_c^1(\R)$, which proves the Lemma, as $\lambda_c$ is assumed to be positive.
\fin\\

Now to prove Proposition \ref{enerpos}, we want to use the Mountain Pass Theorem, so we need to prove that our energy functional satisfies the Palais-Smale Condition.

\begin{Lemma}\label{Palaissmale}
If $(u_n)_n$ is a sequence in $H_c^1(\R)$ such that $E_c[u_n]\leq C$ for all $n\in\N$ and $dE_c [u_n]\to 0$ as $n\to +\infty$ in $(H_c^1)^*$, in the sense that $\lVert dE_c[u_n]\rVert_{(H_c^1)*}\to 0$ as $n\to+\infty$, then there exists a subsequence, that we still call $(u_n)_n$, which converges strongly in $H_c^1(\R)$ toward a solution $u$ of $dE_c[u]=0$.
\end{Lemma}

Here, $(H_c^1)^*$ denotes the dual of the space $H^1_c(\R)$ for the extension of the $L^2_c(\R)$ scalar product. 

\noindent {\bf Proof of Lemma \ref{Palaissmale}}:
As $E_c[u_n]\leq C\textrm{ for all } n\in\N$ and using Lemma \ref{Eminore}, we have
\begin{equation*}
\lVert u_n\rVert_{H_c^1}^2\leq \frac{C+C_1}{\min\{1,\delta\}},
\end{equation*}
which implies that, up to a subsequence, $(u_n)$ converges weakly to $u\in H^1_c(\R)$. Moreover for all $w\in H_c^1(\R)$, $dE_c[u_n](w)\to 0$ as $n\to+\infty$, so
\begin{align*}
0&=\underset{n\to+\infty}{\lim} dE_c[u_n](w)\\
&=\underset{n\to+\infty}{\lim}\int_\R e^{cz}\left\{(u_n)_z  w_z-f(z,u_n)w\right\}dz\\
&= \int_\R e^{cz}\left\{u_z w_z-f(z,u)w\right\}dz.
\end{align*}
Hence $dE_c[u]\equiv 0$.
 
\noindent Now let us prove that $(u_n)$ converges strongly to $u$ in $H^1_c(\R)$ as $n\to+\infty$. We just need to prove that $\ds{\lVert u_n\rVert_{H^1_c(\R)}\to \lVert u\rVert_{H^1_c(\R)}}$ as $n\to+\infty$, since $H_c^1(\R)$ is a Hilbert space. Taking $w=u_n$ we get 
\be\label{weaksolintvn}
\int_\R e^{cz}\left\{(u_n)_z^2-f(z,u_n)u_n\right\}dz=\langle dE_c[u_n], u_n\rangle_{(H^1_c)^*, H^1_c}
\ee
And $\langle dE_c[u_n], u_n\rangle_{(H^1_c)^*, H^1_c}\leq \lVert dE_c[u_n]\rVert_{(H^1_c)^*}\lVert u_n\rVert_{H^1_c}=o(1)$, since $(u_n)$ is bounded in $H_c^1(\R)$. 
Hence, $\langle dE_c[u_n],u_n\rangle\to0$ as $n\to+\infty$.\\ 
As $\langle dE_c[u],u\rangle=0$, we have that 
$$\int_\R e^{cz} u_z^2ds=\int_\R e^{cz}f(z,u)udz.$$
Using the same arguments as in Proposition \ref{minreached} we have that for all $\epsilon>0$,
$$\underset{n\to+\infty}{\lim}\int_\R e^{cz}f(z,u_n)u_ndz\leq \int_\R e^{cz}f(z,u)udz +\epsilon.$$
This inequality and \eqref{weaksolintvn} implies that 
$$\lVert u\rVert_{H^1_c(\R)}\leq \underset{n\to+\infty}{\liminf }\lVert u_n \rVert_{H_c^1(\R)}\leq \underset{n\to+\infty}{\limsup}\lVert u_n \rVert_{H_c^1(\R)}\leq \lVert u \rVert_{H_c^1(\R)}+\epsilon,$$
for all $\epsilon>0$. One has proved the Lemma.\fin
\bigskip\\
\textbf{Proof of Proposition \ref{enerpos}}: As assumed in the Proposition $\underset{u\in H_c^1(\R)}{\min} E_c[u]<0$. 
By Proposition \ref{minreached}, this minimum is reached for some $u_\text{min}\in H_c^1$, and one has $u_\text{min}\not\equiv 0$ since $E_c [0]=0$. 
Lemma \ref{pointhaut} yields that for all $r>0$ sufficiently small, one has $\inf_{\lVert u\rVert_{H_c^1 (\R)}=r} E_c [u]>0$.
Choose $r$ small enough such that $\ds{\lVert u_{min}\rVert_{H_c^1}>r}$.
Then using the Mountain Pass Theorem, there exists $\tilde{u}\in H_c^1$ 
such that $dE_c[\tilde{u}]\equiv 0$ and $E_c[\tilde{u}]\geq \gamma$. We have proved Proposition \ref{enerpos}.\fin
\bigskip\\
We want to prove that we can always find non trivial initial conditions $u_0\not\equiv0$ such that $u$ solutions of \eqref{parabolicgeneral} converge to 0 and to a travelling wave solution.
\begin{Prop}\label{0stable}
Assume that \eqref{0stableeq} holds, then 
\begin{itemize}
\item[$\cdot$] there exists $\overline{u_0}\not\equiv 0$, compactly supported, such that the solution $u$ of \eqref{parabolicgeneral} converges to 0 as $t\to+\infty$,
\item[$\cdot$] there exists $\tilde{u_0}\not\equiv0$, compactly supported, such that the solution $u$ of \eqref{parabolicgeneral} converges to a travelling wave solution as $t\to+\infty$.
\end{itemize}
\end{Prop}

\noindent\textbf{Proof}: We noticed in the previous section that $\lambda_c=\lambda_0+\frac{c^2}{4}$  and if $\lambda_0>0$, then $\lambda_c>0$. We know that there exists a positive function $\phi \in W^{2,p}_\text{loc}(\R)$, for any $1\leq p<+\infty$, such that 
$$-\phi''-c\phi'-f_s(z,0)\phi=\lambda_c\phi \text{ in } \R.$$
Let $w(t,z):=\kappa\phi(z)e^{-\delta t}$ for all $t\geq0$, $z\in\R$, $\kappa>0$, $\delta>0$ some constants that we specify later. Then $w$ satisfies the following equation
$$w_t-w_{zz}-cw_z=(f_s(z,0)+\lambda_c-\delta)w.$$
As $\lambda_c>0$, choosing $\ds{\delta=\frac{\lambda_c}{2}}$, there exists $\kappa>0$ small enough such that
$$w_t-w_{zz}-cw_z\geq f(z,w).$$ 
Thus if $\overline{u_0}\leq \kappa\phi$ in $\R$, using the weak parabolic maximum principle we have that for all $t\geq0$, $z\in\R$,
$$u(t,z)\leq\kappa\phi(z)e^{-\delta t},$$
for some constants $\kappa>0$, $\delta>0$ small enough. 
\bigskip\\
Now we prove the second assertion. As $\underset{u\in H_c^1(\R)}{\inf} \:E_c[u]<0$ there exists $\tilde{u_0}\not\equiv0$, compactly supported such that $E_c[\tilde{u_0}]<0$. As $t\mapsto E_c[u]$ is decreasing, its limit is negative. Thus $u$ converges to a solution of \eqref{elliptproblem} with negative energy, i.e $u$ converges to a travelling wave solution.\\
This proves Proposition \ref{0stable}. \fin
\bigskip\\


\section{Examples and discussion}\label{examples}


\subsection{Numerical simulations}\label{numerics}
In this section we illustrate the behaviour of the solution of the parabolic problem considering different types of reaction terms $f$, different values of $\delta$ and $c$. 
We solve numerically the following problem
\be\label{problemnumeric2}\begin{cases}
\partial_t u-\partial_{zz} u-c\partial_z u=f(z,u), &\text{for } t\in[0,T], z\in[0,L],\\
u(0,z)=e^{-(\frac{z-L/2}{l})^2}, &\text{for } z\in(0,L),\\
u(t,0)=0, &\text{for } t\in[0,T],\\
u(t,L)=0, &\text{for } t\in[0,T],
\end{cases}\ee
where 
\be\label{reactfunc}
f(z,u)=\begin{cases} f_0(u) &\text{if }\frac{L}{2}-\frac{l}{2}<z<\frac{L}{2}+\frac{l}{2},\\
				-\delta u &\text{otherwise},
				\end{cases}
\ee
with $L=300$, $T=150$, $l=30$. We compute this problem using FreeFem++ with $\Delta x=10^{-1}$ and $\Delta t=10^{-1}$. As our solution converges to $0$ as $z\to\pm\infty$ and the initial condition in \eqref{problemnumeric2} is approximatively equal to $1.4\times 10^{-11}$ on the boundary of the domain, we approximate our problem \eqref{problembis} by a Dirichlet boundary value problem with $L$ large enough.  

\subsubsection{Existence of a critical speed}

\begin{figure}[!h]
 \begin{center}
 \includegraphics[scale=0.4]{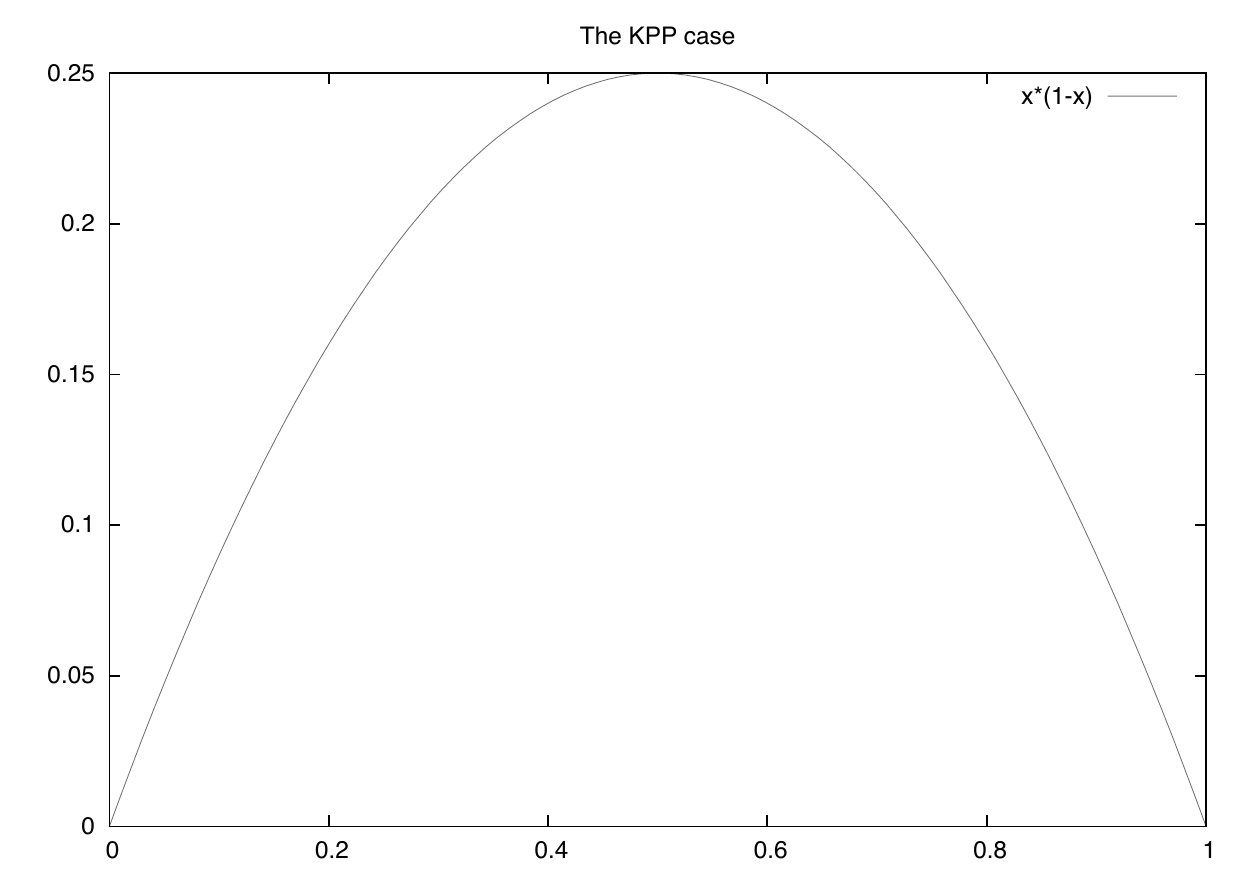}
 \includegraphics[scale=0.4]{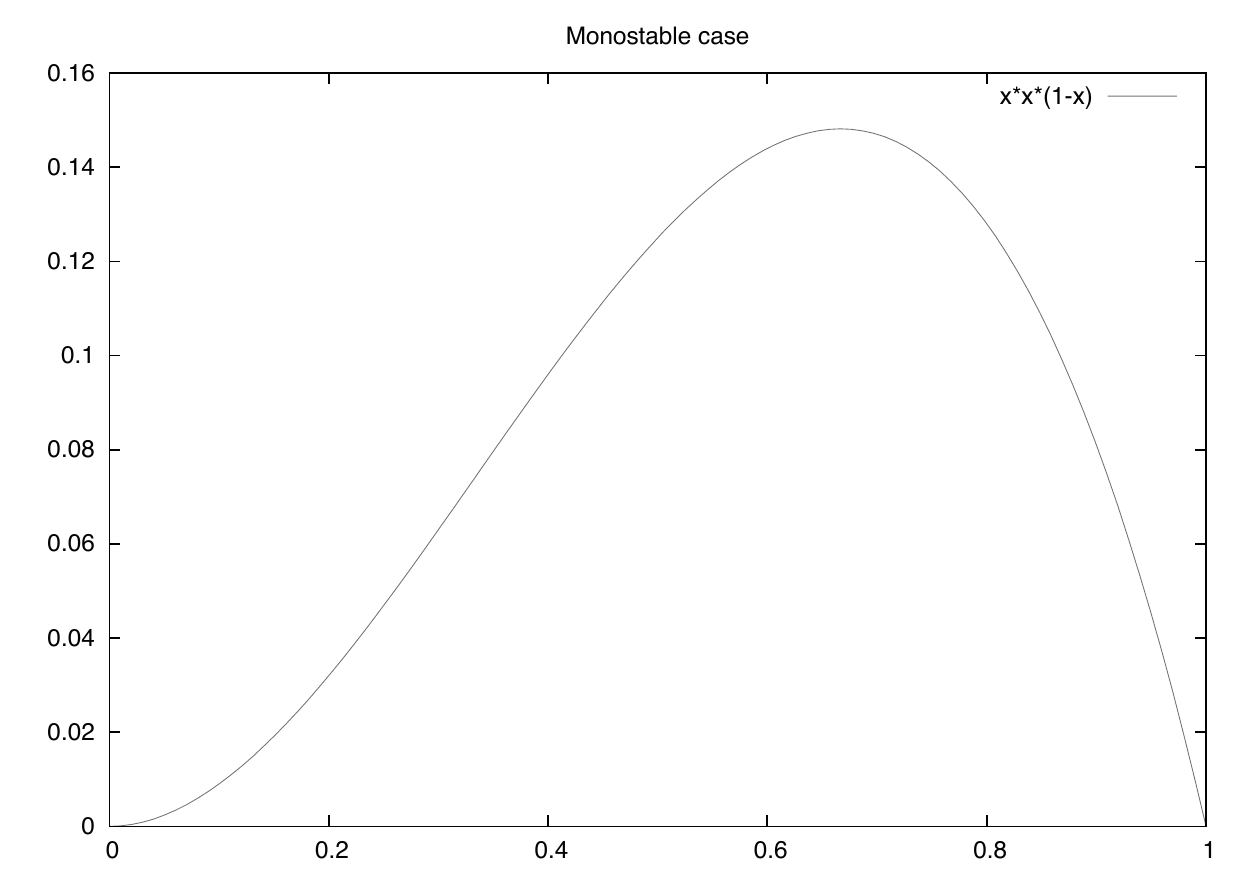}
 \includegraphics[scale=0.4]{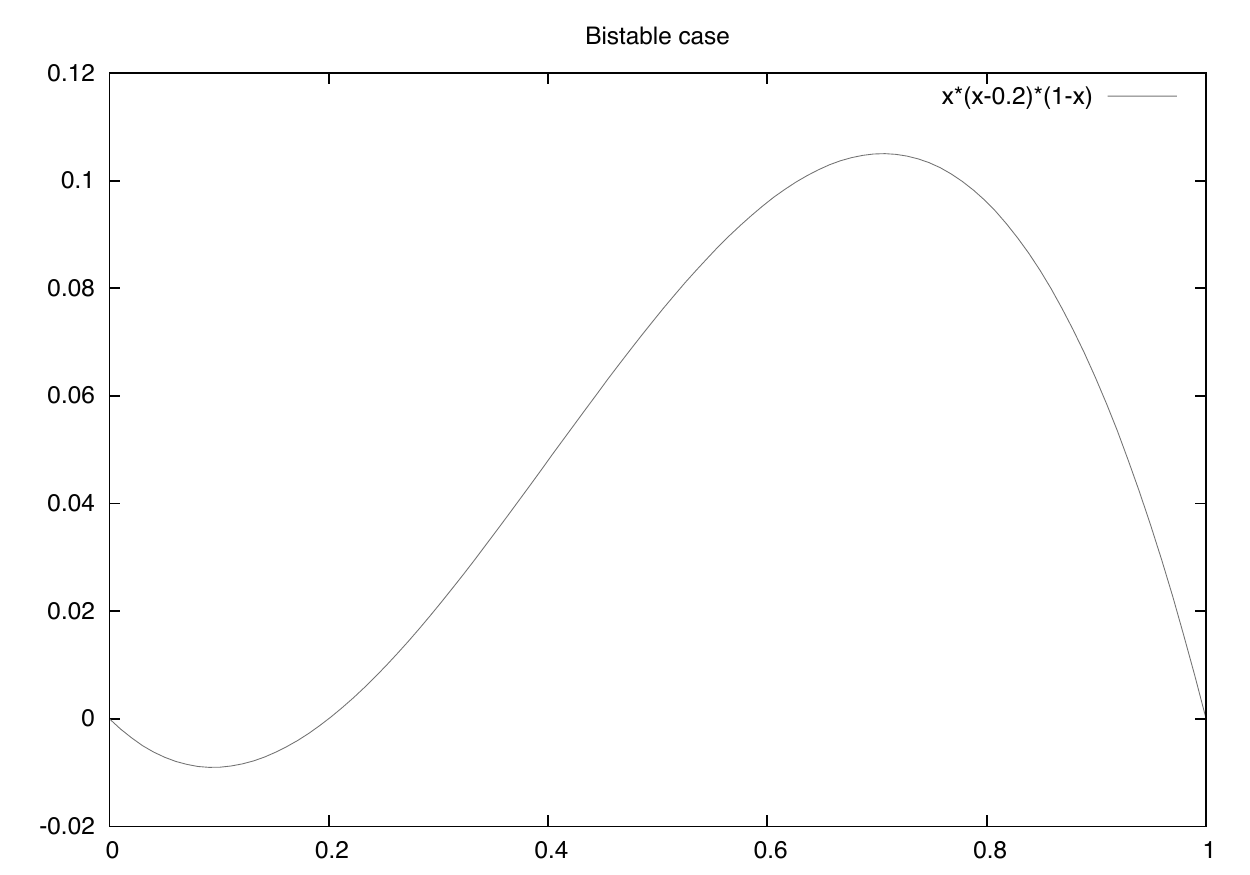}
 \caption{Different types of reaction terms, from left to right: \\
 KPP nonlinearity: $\ds{f_0(u)=u(1-u)}$, Monostable nonlinearity: $\ds{f_0(u)=u^2(1-u)}$ and Bistable nonlinearity: $\ds{f_0(u)=u(1-u)(u-0.2)}$.}\label{reactfuncfig}
 \end{center}
 \end{figure}

We consider three types of reaction function $f_0$: the KPP case, the monostable case and the bistable case (see figure \ref{reactfuncfig}). We restrict our analysis 
to $[0,T]\times[0,L]$ and take $T$ and $L$ large enough to act as if it was $+\infty$.

\noindent In \cite{BDNZ} and \cite{BR1} the authors 
studied the asymptotic behaviour of the parabolic solution and more precisely the existence of non trivial travelling wave solution in the KPP case, 
i.e $\frac{f_0(u)}{u}$ is maximal when $u=0$. The authors proved that there exist travelling wave solutions if and only if $\lambda_0<0$ and $c<2\sqrt{-\lambda_0}$, 
where $\lambda_0$ is the generalized eigenvalue when $c=0$. In other words there exists a critical speed $\ds{c_{lin}=2\sqrt{-\lambda_0}}$ such that 
$\underline{c}=\overline{c}=c_{lin}$ in Theorem \ref{notrivialTW}. In our paper we consider more general nonlinearities $f$ and do not assume that $f$ satisfies the 
KPP property. We proved in Theorem \ref{notrivialTW} that there exists $\underline{c}\leq\overline{c}$ such that there exist travelling  
wave solutions for all $\ds{c<\underline{c}}$ and the only solution of \eqref{elliptproblem} is 0 for all $\ds{c>\overline{c}}$. 
We wonder if in this general framework, there still exists a critical speed, that is, $\ds{\underline{c}=\overline{c}}$. 
We investigate this conjecture numerically in the monostable and bistable case. As the initial data gathers a lot of mass in the favourable area $[\frac{L}{2}-\frac{l}{2};\frac{L}{2}+\frac{l}{2}]$, while it is small in the unfavourable environment, we believe that the solution will converge to a travelling wave solution when it exists for reasonable nonlinearities.  

 \begin{figure}[!h]
 \begin{center}
 \includegraphics[scale=0.2]{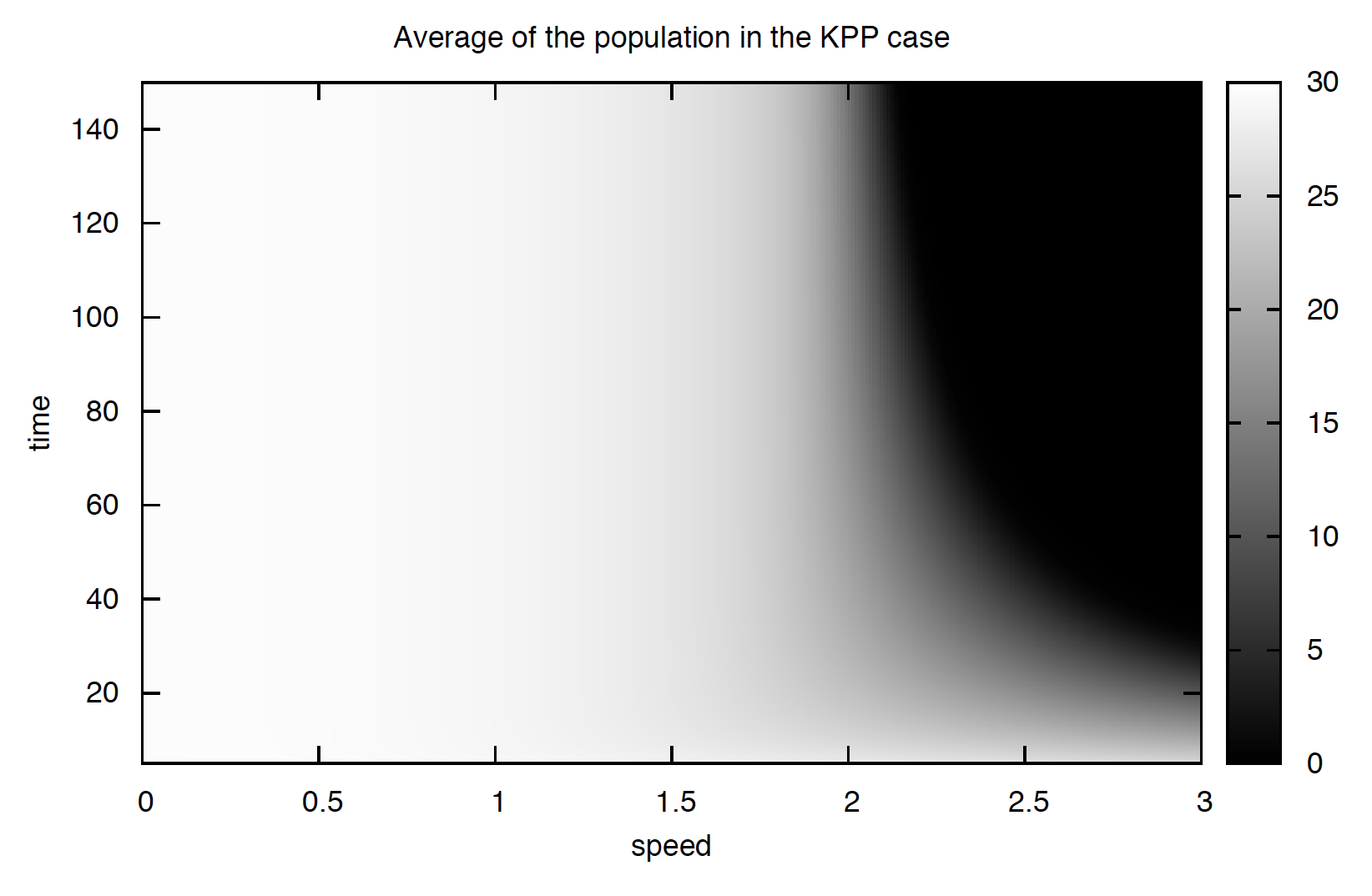}
 \caption{Average of the population $P(t)=\int_0^{L} u(t,x)dx$ for $c\in[0,3]$ in the KPP case (L=300)}\label{critspeedKPPnum}
 \end{center}
 \end{figure}
 \begin{figure}[!h]
 \begin{center}
 \includegraphics[scale=0.2]{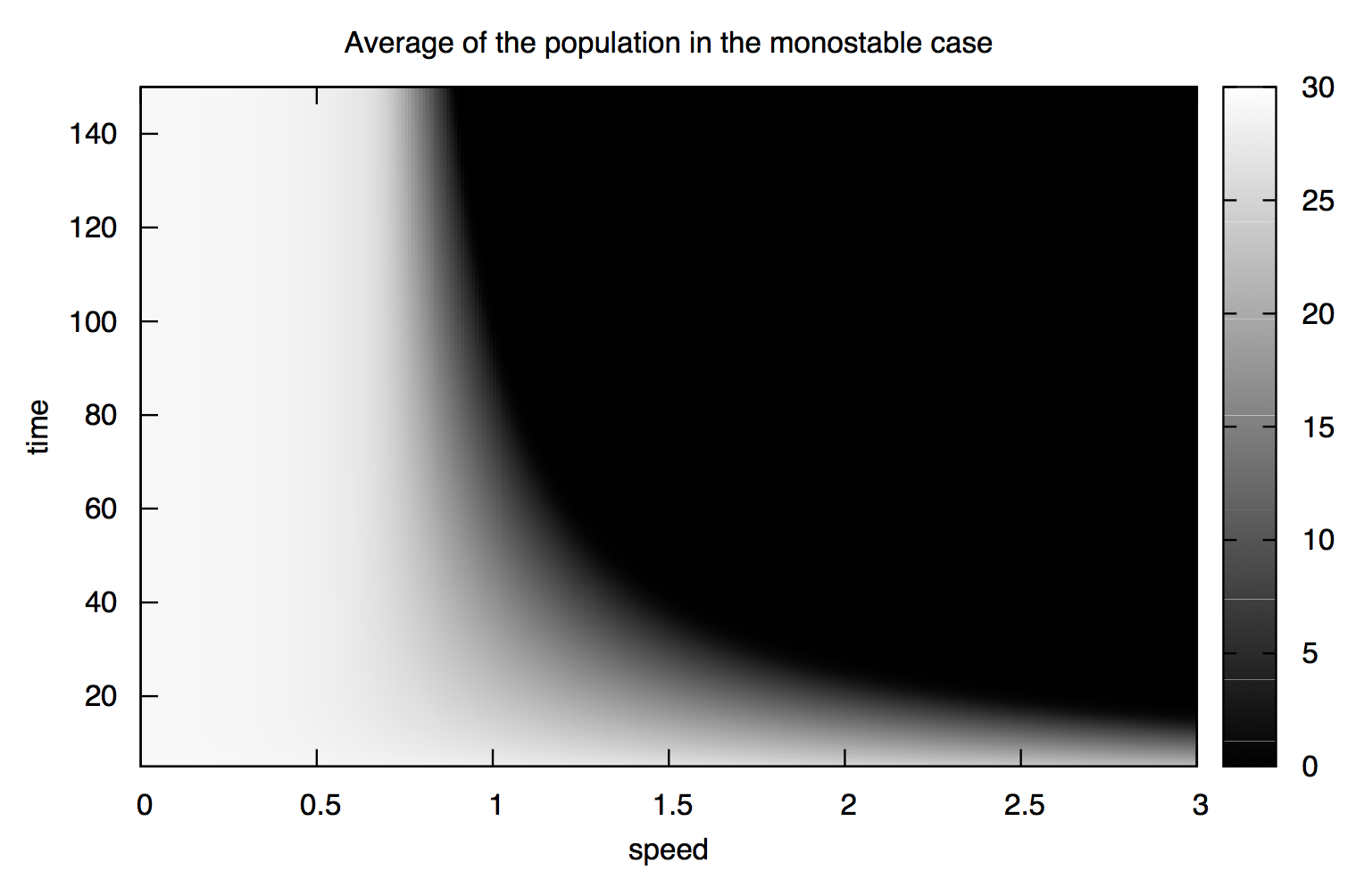}
 \caption{Average of the population $P(t)=\int_0^{L} u(t,x)dx$ for $c\in[0,3]$ in the monostable case (L=300)}\label{critspeedMononum}
 \end{center}
 \end{figure}
 
 \begin{figure}[!h]
 \begin{center}
 \includegraphics[scale=0.2]{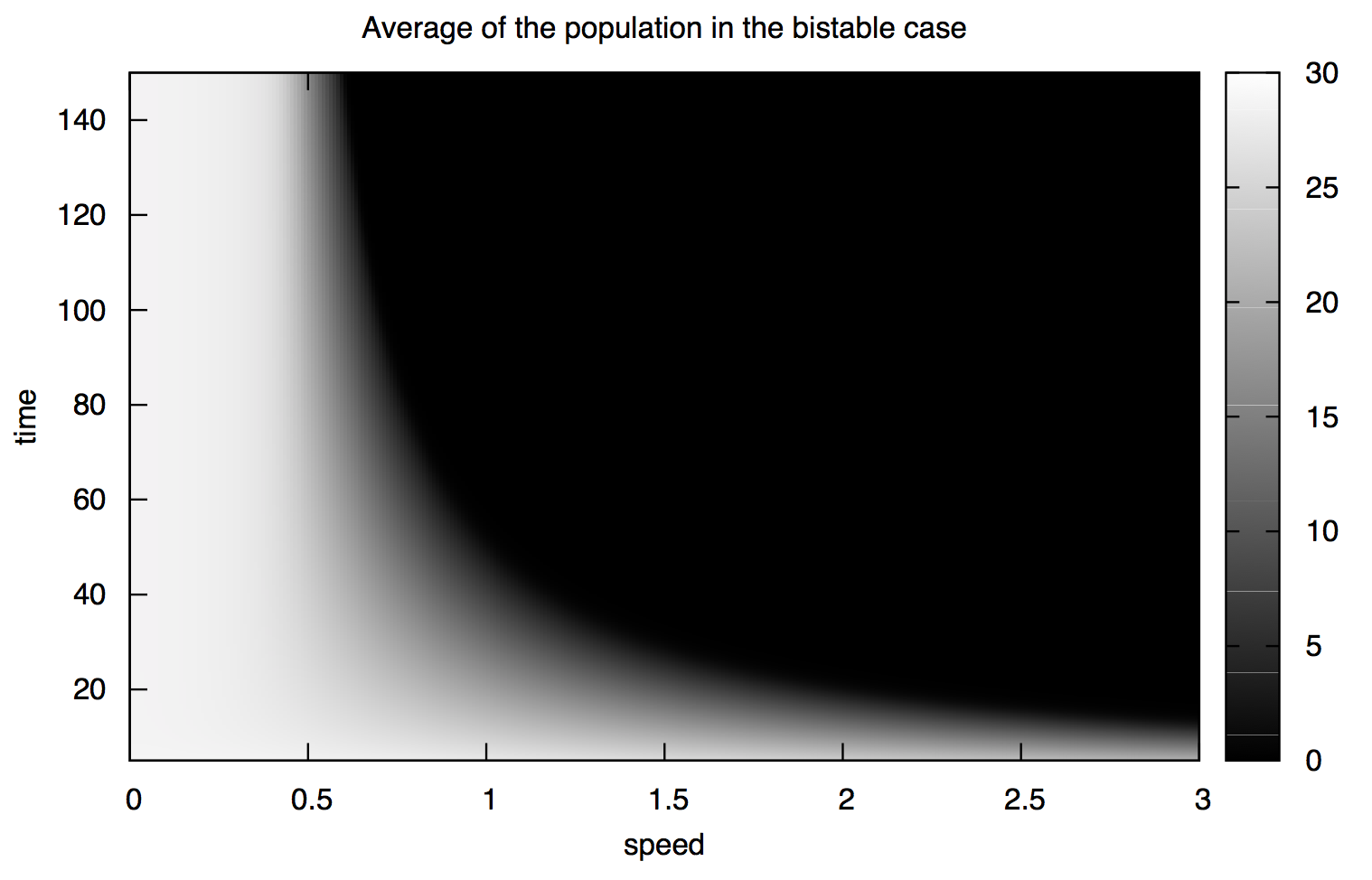}
 \caption{Average of the population $P(t)=\int_0^{L} u(t,x)dx$ for $c\in[0,3]$ in the bistable case (L=300)}\label{critspeedBistablenum}
 \end{center}
 \end{figure}

\noindent The existence of a critical speed has already been introduced in \cite{BRRK}, where the authors highlight some monotonicity of the global population with respect to the speed $c$.\\
Figure \ref{critspeedKPPnum} displays the behaviour proved analytically in \cite{BDNZ, BR1}: there exists a critical speed $\overline{c}$ (around 2) 
such that for $c<\overline{c}$ the population survives whereas for $c>\overline{c}$ the population dies.\\
In Figure \ref{critspeedMononum} and \ref{critspeedBistablenum} one can observe the same phenomenon but for lower critical speeds.
\bigskip\\
Let us also notice that, as proved  in Corollary \ref{critpoint}, we still have existence of travelling waves when $\lambda_c>0$ in the bistable case (Figure \ref{critspeedBistablenum} for $c\in[0,0.4]$).


\subsubsection{Shape of the solution in the moving frame}\label{movingframe}
\begin{figure}[!h]
 \begin{center}
\includegraphics[scale=0.4]{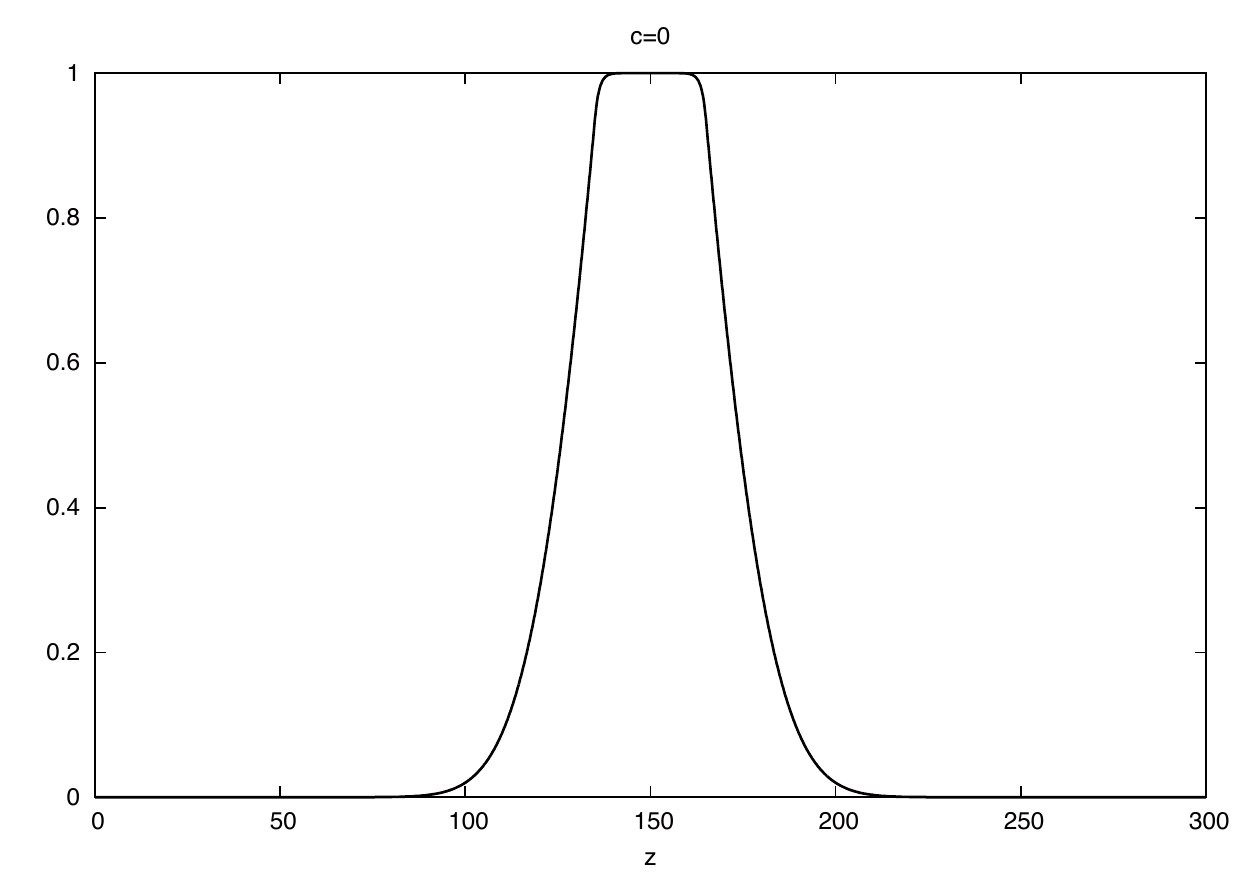}
\includegraphics[scale=0.4]{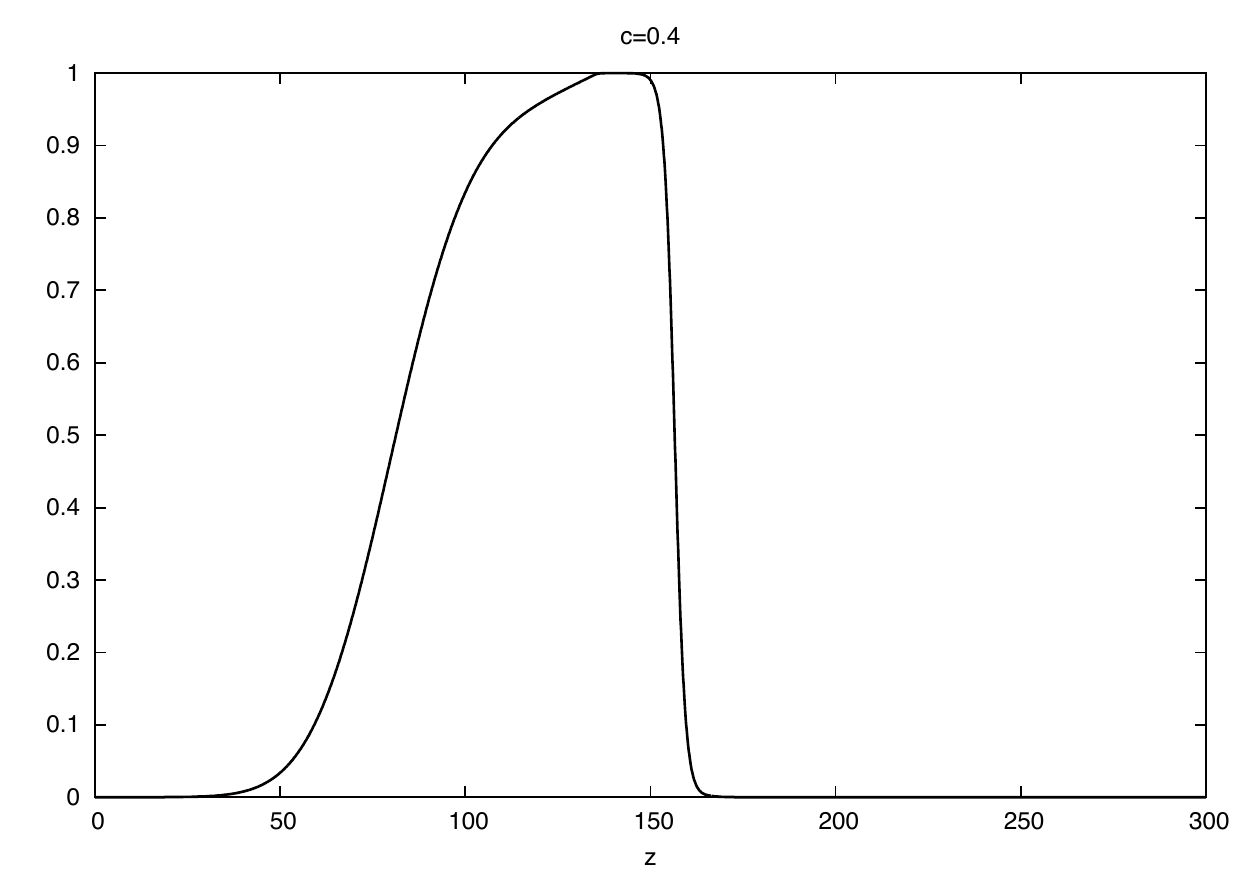}
\includegraphics[scale=0.4]{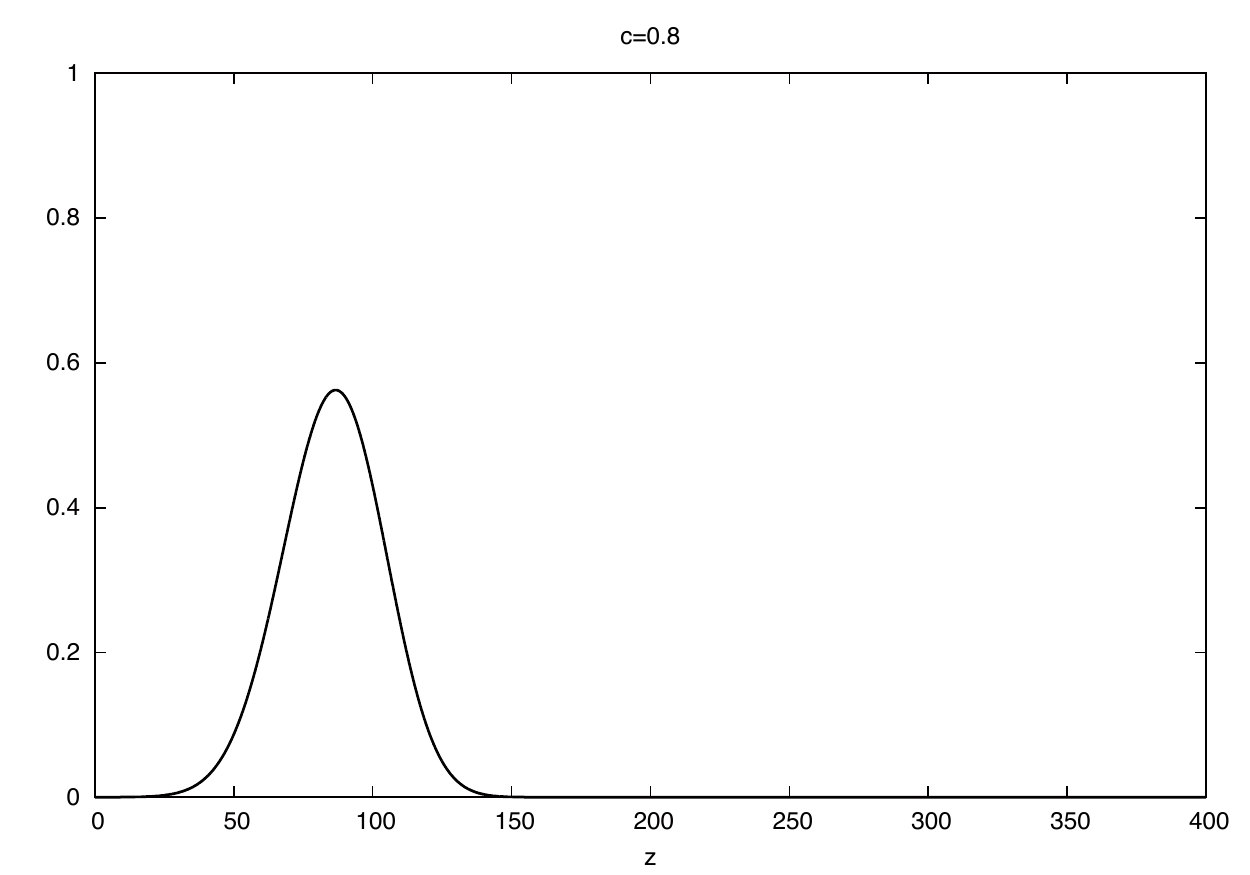}
 \caption{Solution of \eqref{problemnumeric2} for $\delta=0.001$ and $c=0$; $0.4$ and $0.8$ (from left to right) for t=150. For $c=0.8$, the figure on the right displays an intermediate state slowly converging to 0.}\label{r001}
 \end{center}
 \end{figure}

We now investigate the shape of the front when $\delta$ varies and $f$ is bistable, i.e $\ds{f_0(u)=u(1-u)(u-0.2)}$. \\
When $\delta$ is small (figure \ref{r001}), a tail grows at the bottom of the front whereas the transition at the front edge of the front stays sharp when the speed $c>0$ is small enough for the population to survive, as it was already observed by Berestycki et al \cite{BDNZ} for KPP nonlinearity. This tail is created by the movement of the favourable environment, indeed the death rate $\delta$ is too small to kill the population which reproduced quickly in the favourable zone. When the speed is too large the population can not keep tracks with its favourable environment and slowly converges to 0. On the other hand when c=0, both edges of the front become less and less sharp.

\begin{figure}[!h]
 \begin{center}
 \includegraphics[scale=0.5]{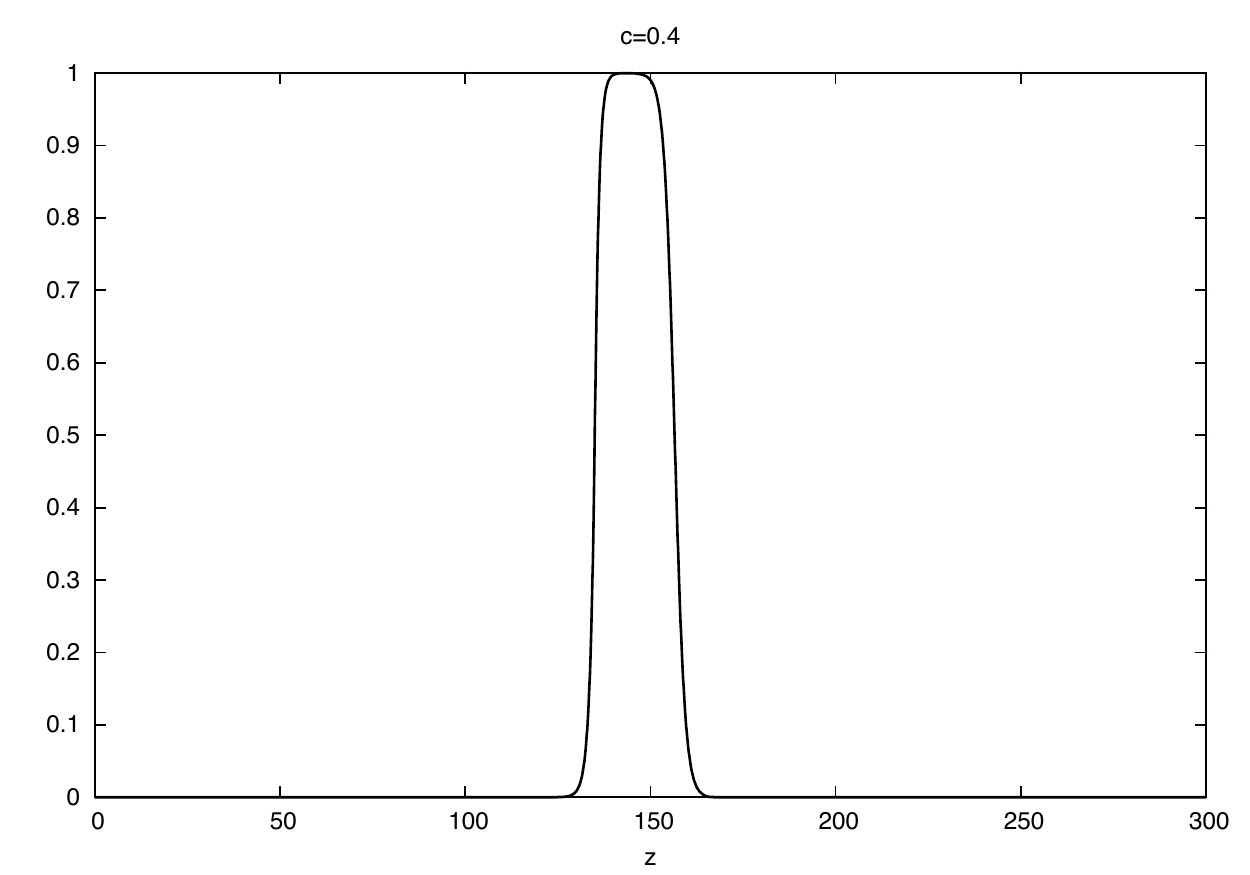}
 \caption{Solution of \eqref{problemnumeric2} for $\delta=1$ and $c=0.4$ for t=150.}\label{r1c04}
 \end{center}
 \end{figure}

\begin{figure}[!h]
 \begin{center}
 \includegraphics[scale=0.5]{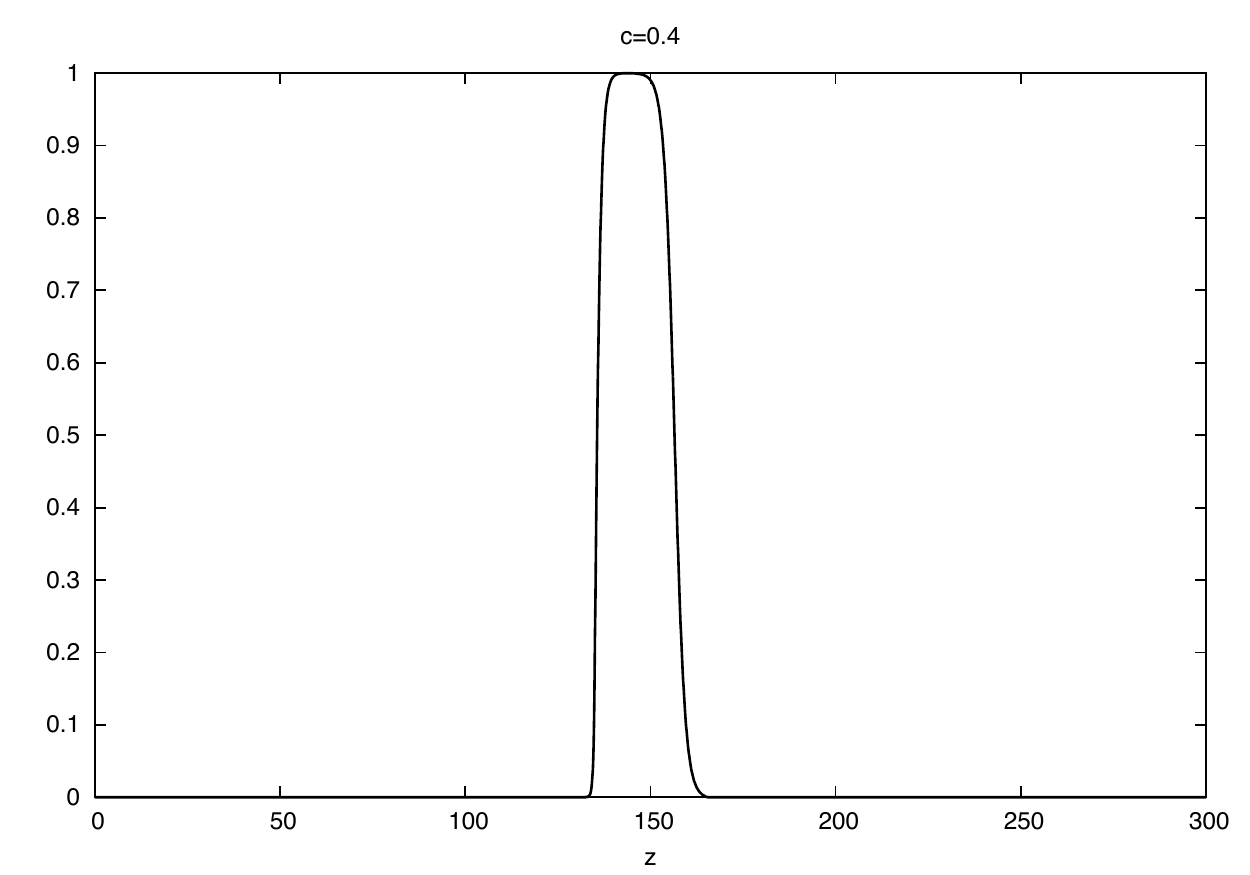}
 \caption{Solution of \eqref{problemnumeric2} for $\delta=10$ and $c=0.4$ for t=150.}\label{r10c04}
 \end{center}
 \end{figure}

\noindent Then we see that when $c>0$ (small enough for the population to survive), both edges of the front become sharper and sharper as $\delta$ increases (Figures \ref{r001}, \ref{r1c04} and \ref{r10c04}).
 
\subsection{Non uniqueness of stable travelling waves}\label{nonuniqprofsection}

We can also build $f$ such that \eqref{elliptproblem} has more than one stable solution with negative energies in the sense that the solutions are local minimisers of the energy functional.

\begin{Prop}\label{energyneg}
There exists $f(z,u)$ satisfying assumptions \eqref{f0}-\eqref{hypenergy}, such that there exist $u^*$ and $v^*$ solutions of \eqref{elliptproblem}  local minimisers of the energy functional with $E_c[v^*]<E_c[u^*]<0$.
\end{Prop}

\begin{figure}[!h]
 \begin{center}
 \includegraphics[scale=0.9]{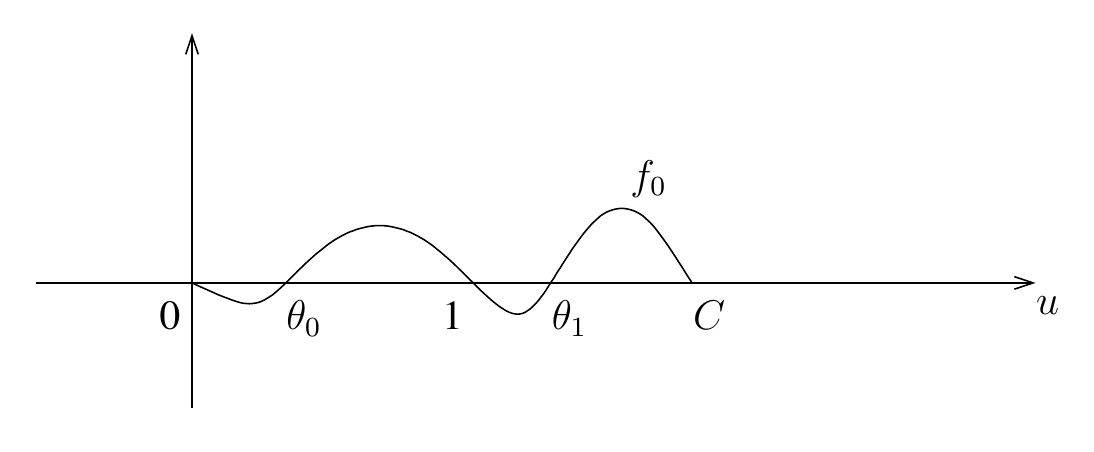}
 \caption{\footnotesize{$f_0$ a multistable function such that $\int_0^1f_0(s)ds>0$, there exist $\ds{0<\theta_0<1<\theta_1<C}$ such that $f(0)=f(\theta_0)=f(1)=f(\theta_1)=f(C)=0$, $\ds{f(s)<0},\quad \text{for }\ds{s\in(0,\theta_0)\cup(1,\theta_1)}$ and $\ds{f(s)>0} \quad \text{for }\ds{s\in(\theta_0,1)\cup(\theta_1,C)}$ with $\int_0^{C}f_0(s)ds>\int_0^1f_0(s)ds$.}}\label{f0multistable}
 \end{center}
 \end{figure}

\noindent Let $f$ be as follow
\be\label{fenergyneg}
f(z,u)=\begin{cases} f_0(u) &\text{if } |z|<R,\\
				-\delta u &\text{otherwise},
	\end{cases}
\ee
where $f_0$ is a multistable function, i.e there exist $0<\theta_0<1<\theta_1<C$ such that 
\begin{gather*}
f(0)=f(\theta_0)=f(1)=f(\theta_1)=f(C)=0,\\
f(s)<0,\quad \text{for }s\in(0,\theta_0)\cup(1,\theta_1),\\
f(s)>0 \quad \text{for }s\in(\theta_0,1)\cup(\theta_1,C),
\end{gather*}
$\int_0^1f_0(s)ds>0$ and $\int_0^{C}f_0(s)ds>\int_0^1f_0(s)ds$ (one can look at Figure \ref{f0multistable} for an example of $f_0$), and $\delta>0$. 

\noindent Such multistable functions have already been used in other frameworks in order to construct multiple stable solutions of semilinear problems (see \cite{BrownBudin}  for example).
 \bigskip\\
\noindent We start with the proof of the following Lemma.

\begin{Lemma}\label{ustar01}
There exists $u^*\in H_c^1(\R)$ a local minimiser of $E_c[u]$ such that $0<u^*<1$ in $\R$, $E_c[u^*]<0$ and $u^*$ is a solution of \eqref{elliptproblem}.
\end{Lemma}

\begin{figure}[!h]
 \begin{center}
 \includegraphics[scale=0.4]{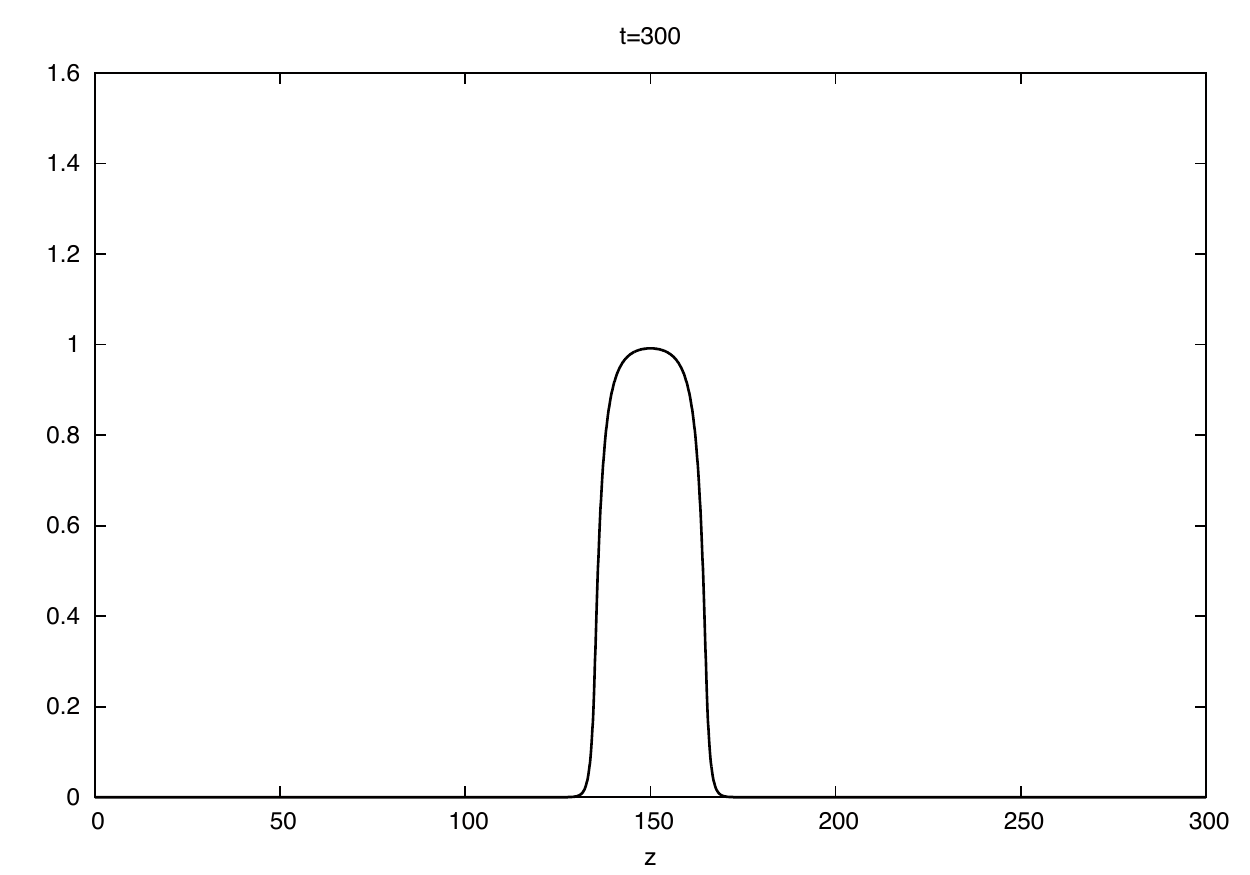}
 \includegraphics[scale=0.4]{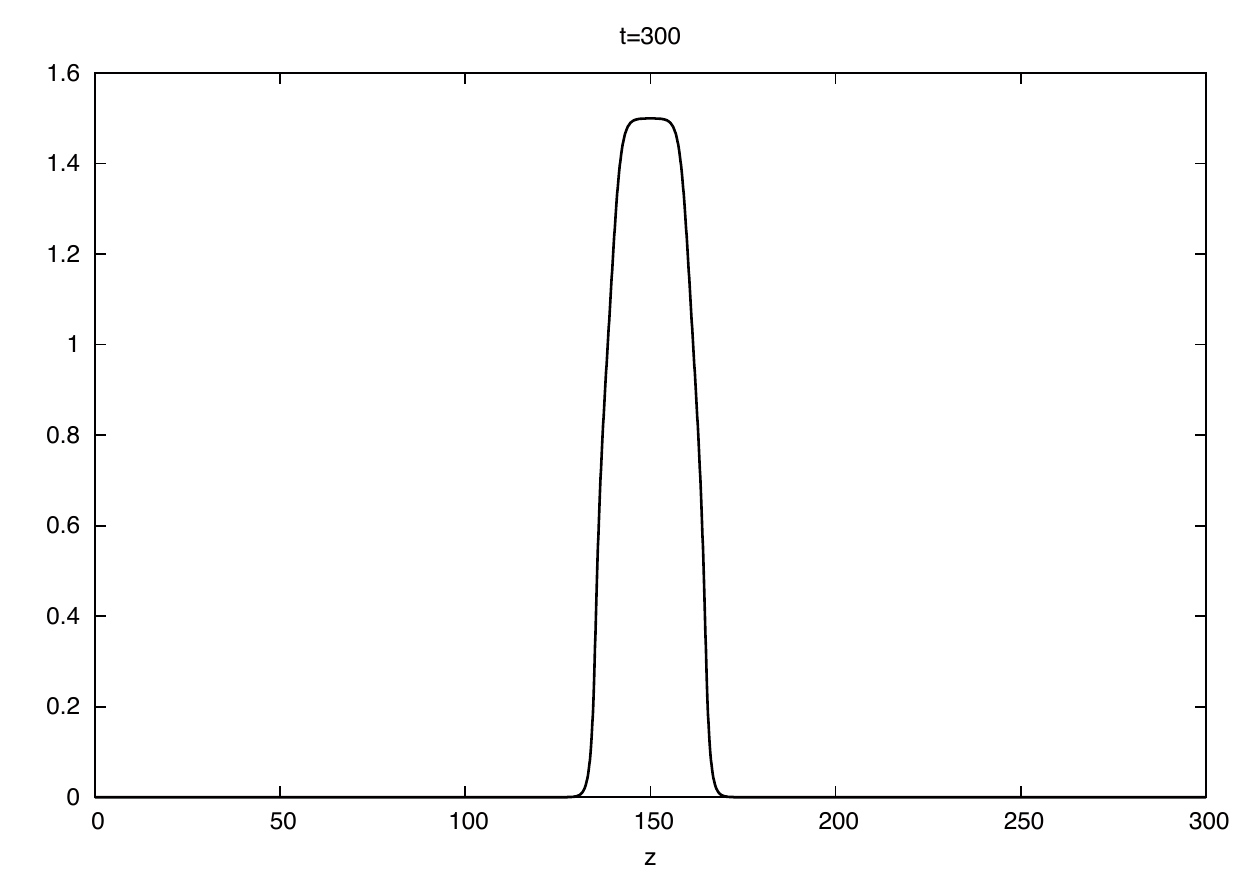}
 \caption{Solution of \eqref{problemnumeric2} for $c=0$ and $t=300$ with $\ds{u_0(x)=e^{-(\frac{z-L/2}{l})^2}}$ for the figure on the left and $\ds{u_0(x)=1.5\times e^{-(\frac{z-L/2}{l})^2}}$ for the figure on the right.}\label{c0}
 \end{center}
 \end{figure}

\begin{figure}[!h]
 \begin{center}
 \includegraphics[scale=0.4]{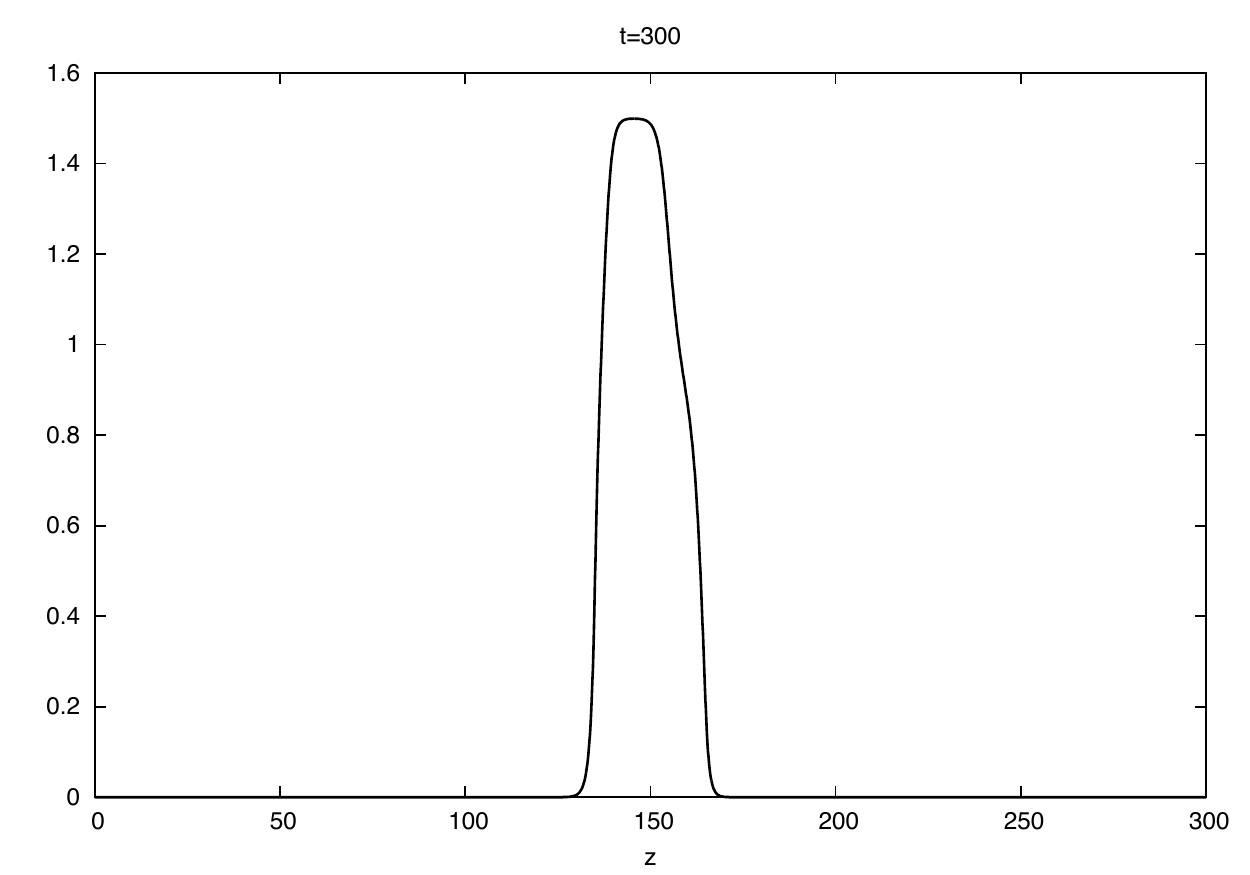}
 \caption{Solution of \eqref{problemnumeric2} for $c=0.2$ and $t=300$ with $\ds{u_0(x)=1.5\times e^{-(\frac{z-L/2}{l})^2}}$.}\label{2c02}
 \end{center}
 \end{figure}

\noindent {\bf Proof of Lemma \ref{ustar01}}: Let us define $f^*$  such that 
\be\label{fstar}
f^*(z,u)=\begin{cases} 0 & \text{if } z\in(-R,R)\text{ and } u\notin [0,1],\\
					f(z,u) &\text{otherwise},
		\end{cases}
\ee
Using Proposition \ref{minreached} we know that there exists $u^*$, travelling travelling wave solution of \eqref{elliptproblem} with $f^*$ for some $c>0$ such that 
$\ds{\underset{u\in H_c^1}{\min}E^*_c[u]=E_c^*[u^*]}$, where $E_c^*$ is the energy functional associated with $f^*$.\\
We know that $u^*\leq 1$ in $\R$ by Remark \ref{Ubdd}. 
Thus $u^*$ satisfies the following equation
$$-(u^*)''-c(u^*)'=f(z,u^*),$$
and 
$$E_c[u^*]=\ds{\underset{u\in H_c^1}{\min}E^*_c[u]}.$$
Taking 
\be
u_{min}(z)=\begin{cases}1 &\text{for all } |z|<R,\\
					0 &\text{for all } |z|>R+1,
				\end{cases}
\ee
such that $u_{min}\in H^1_c(\R)$, one can check that for $R$ large enough $\ds{E_c^*[u_{min}]<0}$, which implies that $\ds{E_c[u^*]=\underset{u\in H_c^1}{\min}E^*_c[u]<0}$. We have proved that there exists a solution $u^*\in H_c^1(\R)$ of \eqref{elliptproblem}, such that $0<u^*$ in $\R$ and $\ds{E_c[u^*]<0}$. Now let us prove that $u^*$ is a local minimiser. Using classical Sobolev injections, there exists $\rho>0$ small enough, such that 
$$\lVert u-u^*\rVert_{H_c^1(\R)}<\rho \quad \implies \quad \lVert u-u^*\rVert_{L^\infty(-R,R)}\leq \theta_1-1.$$
Now let us prove that as soon as $\ds{\lVert u-u^*\rVert_{H_c^1(\R)}<\rho}$, then $\ds{E_c[u]\geq E_c[u^*]}$.
\begin{align*}
E_c[u]&= \int_\R e^{cz}\left\{\frac{(u')^2}{2}-F(z,u)\right\}dz,\\
		&=E_c^*[u]+\int_{-R}^R e^{cz}\left\{F^*(z,u)-F(z,u)\right\}dz.
\end{align*}
As $\lVert u-u^*\rVert_{L^\infty(-R,R)}\leq \theta_1-1$, $f^*(z,u)\geq f(z,u)$ for all $z\in(-R,R)$, thus 
$$\int_{-R}^R e^{cz}\left\{F^*(z,u)-F(z,u)\right\}dz\geq 0.$$
We have proved the Lemma. \fin
\bigskip\\
{\bf Proof of Propostion \ref{energyneg}}: 
Now let us prove that there exists $v^*\in H_c^1(\R)$ solution of \eqref{elliptproblem} such that $E_c[v^*]<E_c[u^*]<0$. Let $u_3$ be as follow,
\be\label{u3}
u_3(z)=\begin{cases} C &\text{if } |z|<R,\\
					0 &\text{if } |z|>R+\epsilon,
	\end{cases}
\ee
such that $u_3\in H_c^1(\R)$. Then 
$$E_c[u_3]=-\left(\int_0^C f_0(s)ds\right)\frac{e^{cR}-e^{-cR}}{c}+\int_{R<|z|<R+\epsilon}\left\{\frac{(u_3'(z))^2}{2}+\frac{\delta u_3(z)^2}{2}\right\}e^{cz}dz.$$
Thus choosing $C$ close enough to 1 and $f_0>>0$ in $(\theta_1+\eta,C-\eta)$ for some $\eta>0$, small, we have 
$$E_c[u_3]<E_c[u^*].$$
Using Proposition \ref{minreached}, we know that there exists $v^*\in H_c^1(\R)$ such that 
$$\ds{E_c[v^*]=\underset{u\in H_c^1(\R)}{\min} E_c[u]}\leq E_c[u_3].$$ One has proved Proposition \ref{energyneg}.\fin
\bigskip\\
We now illustrate the previous results. Choosing a specific reaction term
$$\ds{f_0(u)=u(1-u)(u-0.2)(1.1-u)(1.5-u)}$$ 
and an appropriate initial condition we get different convergence results as one can see in Figures \ref{c0} and \ref{2c02}. We computed the same problem \eqref{problemnumeric2} that in section \ref{movingframe}, with $\delta=1$ and $\ds{f_0(u)=u(1-u)(u-0.2)(1.1-u)(1.5-u)}$.
In the first figure (Figure \ref{c0}), one can see that depending on the initial condition, we get two different fronts but with a similar shape with sharp edge on both sides. On the other hand when $c>0$ the front edge takes the shape of a stairs, indeed in the favourable environment the population moves rapidly to 1 but need more time to grow from 1 to 1.5.

\clearpage

%
\vskip 20pt
\bibliography{bibliothese}
\bibliographystyle{plain}
\end{document}